\documentclass[twoside,12pt]{article}
\usepackage{latexsym}
\usepackage{amsmath}
\usepackage{amssymb}
\usepackage{graphics}
\usepackage{lscape}
\usepackage{hhline}

\newcommand\hiii{\mathbf H^3}
\newcommand\hiv{\mathbf H^4}
\newcommand\hn{\mathbf H^n}

\newcommand\en{\mathbf E^n}
\newcommand\sn{\mathbf S^n}
\newcommand\siii{\mathbf S^3}
\newcommand\eiii{\mathbf E^3}

\newcommand\rii{\mathbf R^2}
\newcommand\riii{\mathbf R^3}
\newcommand\riv{\mathbf R^4}

\newcommand\mbar{\overline M}

\DeclareMathOperator{\Isom}{Isom}

\DeclareMathOperator{\intr}{int}

\newcommand{\arrtop}[2]{\xrightarrow[#2]{#1}}

\newcommand\bd{\partial}
\newcommand\qed{\Box}

\title{On converting a side-pairing to a handle decomposition}
\author{Dubravko Ivan\v si\'c}
\date{Murray State University}

\begin{document}

\newtheorem{theorem}{Theorem}[section]
\newtheorem{proposition}[theorem]{Proposition}
\newtheorem{remark}[theorem]{Remark}
\newtheorem{example}[theorem]{Example}

\maketitle

\begin{abstract}
We give a method for obtaining a handle decomposition of an $n$-manifold if the manifold is given
by isometric side-pairings of a polyhedron in $\en$, $\sn$ or $\hn$.
Every cycle of $k$-faces on the polyhedron corresponds to an $(n-k)$-handle of the manifold.

Two applications of the method are given.
One helps recognize when a noncompact hyperbolic 3-manifold is a complement of a link in $S^3$ (and
automatically produces the link diagram), the other shows that a topological $S^4$ described by the author
in \cite{Ivansic3} is diffeomorphic to the standard differentiable $S^4$.
\vskip8pt
\noindent
{\it MSC:} 57M50, 57M25, 57Q45
\vskip8pt
\noindent
{\it Keywords:} handle decomposition, hyperbolic manifold, link complement,
simply-connected closed 4-manifold
\end{abstract}

\section{Introduction}
\label{introduction}

Let $X$ be any of the three constant-curvature spaces $\en$, $\sn$ or $\hn$, and let $G$ be a
discrete subgroup of isometries of $X$.
By a geometric manifold we mean a manifold of the form $M=X/G$.

Many examples of geometric manifolds are given through side-pairings of a polyhedron
$P\subset X$, this being a convenient and topologically revealing way of describing a manifold.
On the other hand, general manifolds are often given using a handle decomposition, which
lends itself to manipulation and simplification through handle moves.

In this paper we give a method that converts a polyhedron-side-pairing representation of a manifold
into a handle decomposition of the manifold.
The method associates every cycle of $k$-faces in the polyhedron to an $(n-k)$-handle in the
handle decomposition.
While the method works in any dimension, it is most interesting to us in dimensions $n=3,4$, where
we give two applications.

In section~\ref{conv3} we motivate and illustrate the method by describing it in dimension~3, 
where it is easily understood.

Section~\ref{ident3} provides an application of the method to hyperbolic 3-manifolds.
Many examples of finite-volume hyperbolic manifolds $M$ are known to be complements of links
in the 3-sphere.
However, proving that a particular manifold is a complement of a particular link is often demanding
and pushes the limits of intuition.
Furthermore, proofs that the author has seen usually require that the link is known before one
executes the proof.
(The only procedure the author is aware of that does not require this is
described in Francis' book \cite{Francis}, however, the procedure is significantly restricted by
the type of side-pairings it works for.)
We use the method of Section~\ref{conv3} to obtain a handle decomposition of a given hyperbolic manifold.
Using handle moves one can easily show that the manifold is a complement of a link
in the 3-sphere, while the handle moves produce the diagram of the link as the computation
progresses.
This procedure has worked in a straightforward way on all the standard examples
(complements of the figure-8 knot, the Whitehead link and the Borromean rings) and some
less standard ones, like those in \cite{Wielenberg}.

In Section~\ref{convgen} we justify the conversion method for all dimensions.

Section~\ref{diagram4} details how to get handle decomposition diagrams in dimension~4.

Section~\ref{ident4} gives an application of the conversion method in dimension~4.
J.~Ratcliffe,  S.~Tschantz and the author have found a dozen examples (see \cite{Ivansic3, Ivansic4})
of noncompact hyperbolic 4-manifolds that are complements of varying numbers of tori and Klein
bottles in a topological 4-sphere $N$.
We work out the handle decomposition of one $N$ in order to show that it is diffeomorphic to the standard
differentiable 4-sphere, which the original proof was not equipped to do.

As a matter of fact, the author's motivation for developing the conversion method was the problem of
whether the topological 4-spheres found in \cite{Ivansic3, Ivansic4} were diffeomorphic to the
standard 4-sphere.
The dimension-3 application from Section~\ref{ident3} was found afterwards.

\begin{figure}
\begin{center}
\resizebox{2.5in}{!}{\includegraphics{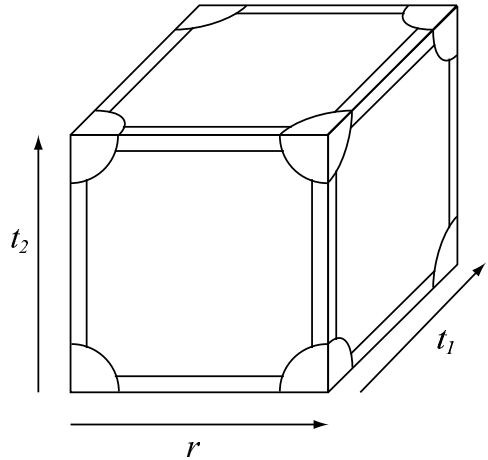}}
\caption{Cube with side-pairing, neighborhoods of faces}
\label{faceneighborhoods}
\end{center}
\end{figure}

\vfill

\section{Conversion in dimension 3}
\label{conv3}

Let $P$ be a polyhedron in $X=\hiii$, $\eiii$ or $\siii$ with a side-pairing defined on it that gives a
geometric manifold $M$.
In Fig.~\ref{faceneighborhoods} a cube is drawn as an example: its top and bottom and front and back sides
are paired by a translation, while the left and right sides are paired by a translation followed by a $180^\circ$
rotation around the translation vector.

Select neighborhoods (for example, $\epsilon$-neighborhoods) around vertices and edges like in
Fig.~\ref{faceneighborhoods}.
The neighborhoods should match via the side pairing.
Let $V_1,\dots, V_m$ be neighborhoods of a cycle of vertices $\{v_1,\dots,v_k\}$ (a cycle of faces
comprises all the faces of $P$ that are identified by the side-pairing).
Then $V_1\cup\dots\cup V_m$ assembles into a ball $V$ in $M$.
In our example, all the vertices are in the same cycle, and $V_i$ is an eighth of a ball.
Eight such pieces, of course, assemble in a ball.

Removing neighborhoods of all vertices from $P$ removes parts of the neighborhoods of the edges.
Let $E_1,\dots,E_n$ be the truncated neighborhoods of a cycle of edges $e_1,\dots,e_n$.
Then $E_1\cup\dots\cup E_n$ assembles into a solid cylinder around a truncated edge, which can also
be viewed as a 3-ball $E$ in $M$.

\begin{figure}
\begin{center}
\resizebox{3.5in}{!}{\includegraphics{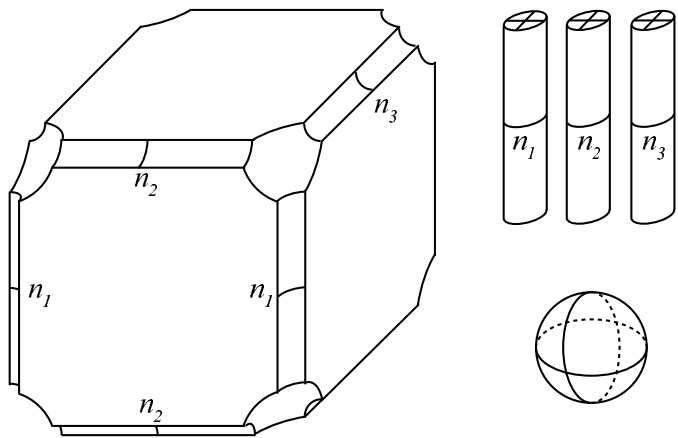}}
\caption{Handles as assemblies of face neighborhoods}
\label{truncatedcube}
\end{center}
\end{figure}

Let $H_1$ be the solid obtained by removing neighborhoods of  vertices and truncated
neighborhoods of edges from $P$.
On the surface of $H_1$ it is the truncated sides that get identified, representing pairwise-identified
disjoint disks, so $H_1$ projects to a handlebody $H$ in $M$ under the quotient map $P\to M$.
The feet of the 1-handles of $H$ are the truncated sides on $H_1$ (see \cite{Gompf-Stipsicz} for basics
of handles and handle decompositions).
Now, the ball $E=D^2\times D^1$ from above is attached to $H$ along
$\partial D^2\times D^1$, making it a 2-handle of $M$.
In our example, there are three cycles of edges, and the visible portions of attaching circles
$n_1$, $n_2$ and $n_3$ of the corresponding 2-handles are shown in Fig.~\ref{truncatedcube}.
Of course, the ball $V$ from above may be viewed as $V=D^3\times D^0$, and it attaches to 
the 0-, 1- and 2-handles along $\partial D^3\times D^0$, making it a 3-handle.

If $P$ is a polyhedron in $\hiii$ with some ideal vertices, the procedure works the same way,
except, instead of removing a neighborhood of the vertex we remove a horoball centered at
the ideal vertex.

Therefore, to get a handle decomposition diagram (pairs of disks in $\rii$ representing feet
of 1-handles, curves outside of the disks representing attaching circles of 2-handles), 
do the following:

\begin{figure}
\begin{center}
\resizebox{3in}{!}{\includegraphics{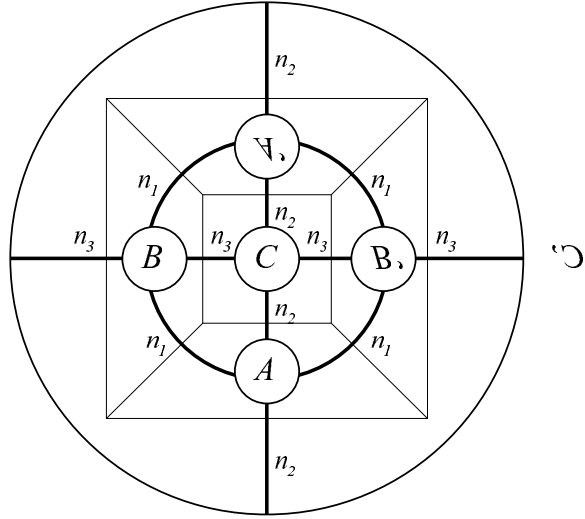}}
\caption{Handle decomposition for the side-pairing from Fig.~1}
\label{convertedcube}
\end{center}
\end{figure}

\begin{itemize}
\item[---]
Project the surface of the polyhedron $P$ to $\rii\cup\infty$ and draw its decomposition into sides.
(If the polyhedron has ideal vertices, one may draw them as empty circles.)
\item[---]
Draw a disk inside every side that represents one of the feet of a 1-handle
(paired sides correspond to feet of 1-handles).
One of the disks may  be the outside of the diagram since a sphere (the surface
of $P$) was projected to $\rii$.
\item[---]
If two sides are adjacent along an edge $e$, draw an arc crossing $e$ once between the disks
corresponding to the sides.
The union of arcs crossing edges that are in the same cycle comprise the attaching circle for a 2-handle.
\item[---]
Attention needs to be paid to how disks (feet of 1-handles) are identified, as the
transformation that identifies them depends on the transformation that identifies the corresponding
sides of $P$.
(We do not assume that the feet of 1-handles are identifed by a reflection in the bisector
of the centers, as is common in handle-decomposition diagrams.)
\item[---]
It is not necessary to keep track of 3-handles, since there is only one way then attach them.
Furthermore, if the polyhedron is hyperbolic and has only ideal vertices, there are no 3-handles.
However, if some of the vertices are real and some ideal, it may be useful to note where on the
diagram  the 3-handles attach.
If necessary, one might put a full circle in $\rii$ wherever there was a real vertex to indicate
that the boundary of a 3-ball is attached to that section of $\rii$, and put an empty circle wherever
there was an ideal vertex to signify that this part of $\rii$ becomes a part of the boundary of the
manifold.
\end{itemize}

Fig.~\ref{convertedcube} illustrates the process above for the cube example at the beginning of the section.
The letters inside the disks suggest the map that pairs the two disks, for example, $A$ and $A'$ are
paired by a reflection in their bisector, while $B$ and $B'$ are paired by a reflection in the bisector,
followed by a rotation by $180^\circ$.

\section{Identifying hyperbolic 3-manifolds as link complements in the 3-sphere}
\label{ident3}

In this section, we apply the conversion method of \S~\ref{conv3} to
illustrate a procedure that attempts to show that a finite-volume noncompact hyperbolic manifold
is the complement of a link in the 3-sphere.
If the procedure is carried out successfully, it also produces the link diagram.

We will use  Wielenberg's example 4 from \cite{Wielenberg}.
In that paper, the following algebraic theorem of Riley's \cite{Riley} is used to determine that  certain
hyperbolic manifolds $M$ are complements of links $S^3-L$: if $\pi_1 M$ is anti-isomorphic to
$\pi_1 (S^3-L)$, then $M\cong S^3-L$.
In order to verify the anti-isomorphism, however, the link has to be known in advance to get the presentation
of $\pi_1 (S^3 - L)$.
Our procedure produces the link diagram as it is carried out.

\begin{figure}
\begin{center}
\includegraphics{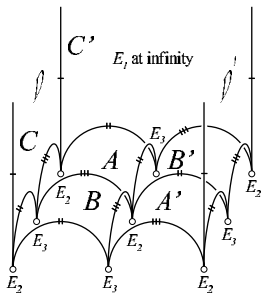}
\caption{Wielenberg's side-pairing on a hyperbolic polyhedron}
\label{wielenbergpolyhedron}
\end{center}
\end{figure}

\begin{figure}
\begin{center}
\includegraphics{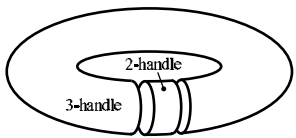}
\caption{Attaching a solid torus along boundary is like adding a 2-handle and a 3-handle}
\label{torusas2handle}
\end{center}
\end{figure}

The hyperbolic manifold comes from pairing the sides of the polyhedron $P$ pictured in the upper
half-space model in Fig.~\ref{wielenbergpolyhedron}.
The vertical sides $C$, $C'$ and $D$, $D'$ are paired by translations.
Sides $A$ and $A'$ are paired by a reflection in the vertical plane passing through the point
where $A$ and $A'$ touch.
Side $B$ is sent to $B'$ by a reflection in the vertical plane that slices $B$ and $B'$
in half, followed by a translation that slides $B$ to $B'$.

A finite-volume orientable noncompact hyperbolic 3-manifold $M$ is diffeomorphic to the interior of a
compact 3-manifold $\mbar$, whose boundary components are all tori.
If solid tori are glued onto the boundary components and the result is a 3-sphere, then $M$
is diffeomorphic to $S^3-L$, where $L$ is the collection of the center circles of the solid tori we added.

\begin{figure}
\begin{center}
\resizebox{4in}{!}{\includegraphics{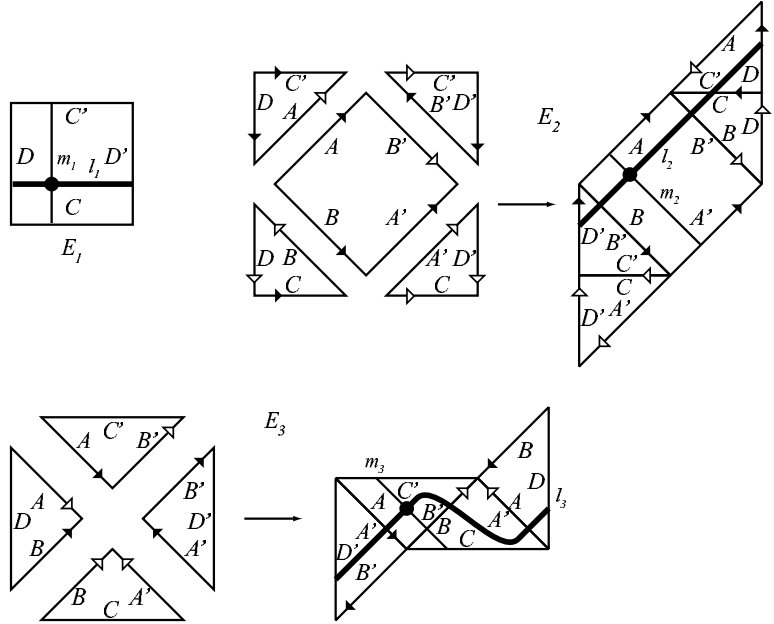}}
\caption{Finding suitable meridians in boundary components}
\label{wielenbergvertexlinks}
\end{center}
\end{figure}

As Fig.~\ref{torusas2handle} suggests, gluing a solid torus to a component $T^2$ of $\bd\mbar$ is the same as
attaching a 2-handle and a 3-handle to $\mbar$.
The attaching circle of the 2-handle can be any nontrivial simple closed curve on $T^2$.
The  components of $\bd\mbar$ are assembled from polygons, called vertex links, that are intersections of small
enough horospheres centered at ideal vertices with the polyhedron $P$.
In our example, the vertex links are $45^\circ$-$45^\circ$-$90^\circ$ triangles and squares.
The three cycles of ideal vertices, $E_1$, $E_2$ and $E_3$ are indicated in Fig.~\ref{wielenbergpolyhedron},
and in Fig.~\ref{wielenbergvertexlinks} the vertex links from each cycle are drawn together and it is shown
how they assemble into parallelograms that give rise to toral boundary components of $\mbar$.

For every boundary component $T^2$ we now choose two curves representing generators of
$\pi_1 T^2$.
One will serve as the attaching circle of the 2-handle, making it a meridian of the attached solid torus.
The other automatically becomes a longitude of the solid torus, thus isotopic to its center circle.
In Fig.~\ref{wielenbergvertexlinks}, the attaching circle is the thinner arc $m_i$ and the longitude is the
thicker arc $l_i$, $i=1,2,3$.

When choosing the attaching circle, choose a curve in $T^2$ as short as possible (in its Euclidean metric).
If the length of the attaching circle is more than $2\pi$, the $2\pi$-theorem on hyperbolic Dehn surgery
(\cite{Bleiler-Hodgson}) asserts we will get a hyperbolic manifold.
Thus, if $\bd\mbar$ has only one component, we will have failed to produce $S^3$.
If $\bd\mbar$ has several components, it is possible that some combination of long and short
attaching circles still produces $S^3$, but chances are probably better the greater the number of
shorter ones are chosen.

Let $M_W$ now denote the manifold resulting from the side-pairing on the polyhedron above.
Since there are three cycles of ideal vertices, $\bd\mbar_W$ will have three components.
Step~0 of Fig.~\ref{handlecancellation1} shows  the handle decomposition of $\mbar_W$,
obtained using the conversion method from~\S\ref{conv3}.
The feet $A$, $A'$, $C$, $C'$ and $D$, $D'$ of 1-handles are all identified by a reflection
in the perpendicular bisector of the line connecting their center.
The feet $B$ and $B'$ are identified by a reflection in the line joining their centers, followed by a translation
that moves $B$ to $B'$, so that the arrows drawn inside match up.
Attaching circles coming from cycles of edges are labeled I, II and III.
The attaching circles that we chose in Fig.~\ref{wielenbergvertexlinks} are also drawn in and labeled
$m_1$, $m_2$, and $m_3$.
Their corresponding longitudes $l_1$, $l_2$ and $l_3$, are drawn as thick curves.

\begin{figure}
\begin{center}
\includegraphics{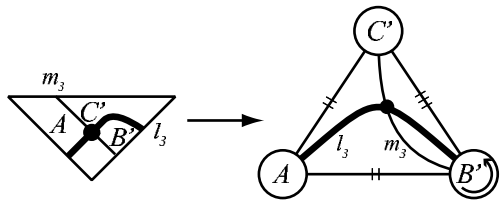}
\caption{Converting one diagram to another}
\label{linkstohandle}
\end{center}
\end{figure}

Fig.~\ref{linkstohandle} shows to make the easy correspondence between a triangle appearing in
Fig.~\ref{wielenbergvertexlinks} and the section of the boundary of the handlebody in step~0
of Fig.~\ref{handlecancellation1} necessary to draw in the longitudes and meridians.
The handle decomposition of $\mbar_W$ does not have any 3-handles, since $P$ did not have
any real vertices.
However, closing off $\bd\mbar_W$ with three solid tori adds three 3-handles.

\begin{figure}
\begin{center}
\resizebox{\textwidth}{!}{\includegraphics{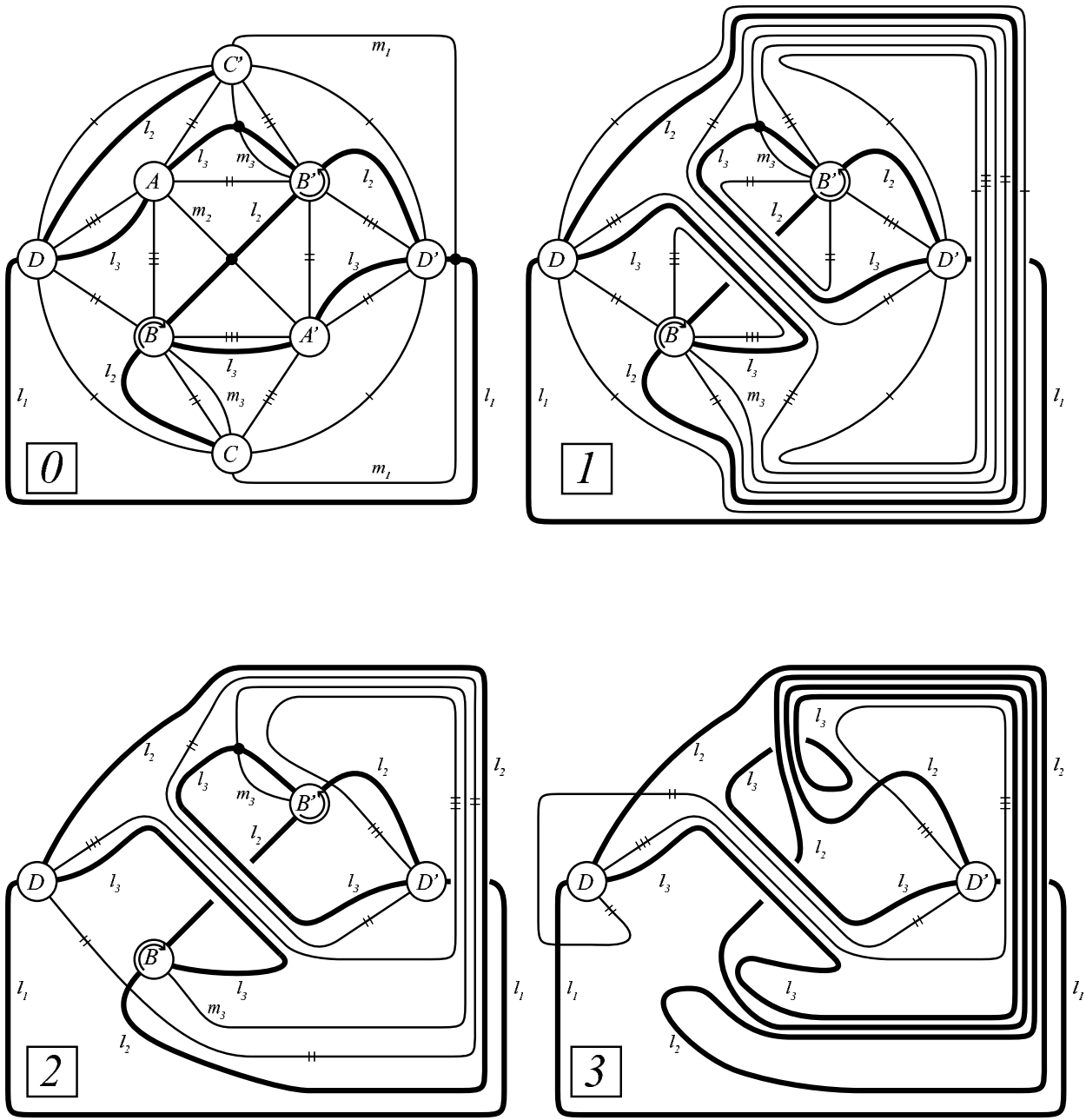}}
\caption{Handle moves, steps 0--3}
\label{handlecancellation1}
\end{center}
\end{figure}

\begin{figure}
\begin{center}
\resizebox{\textwidth}{!}{\includegraphics{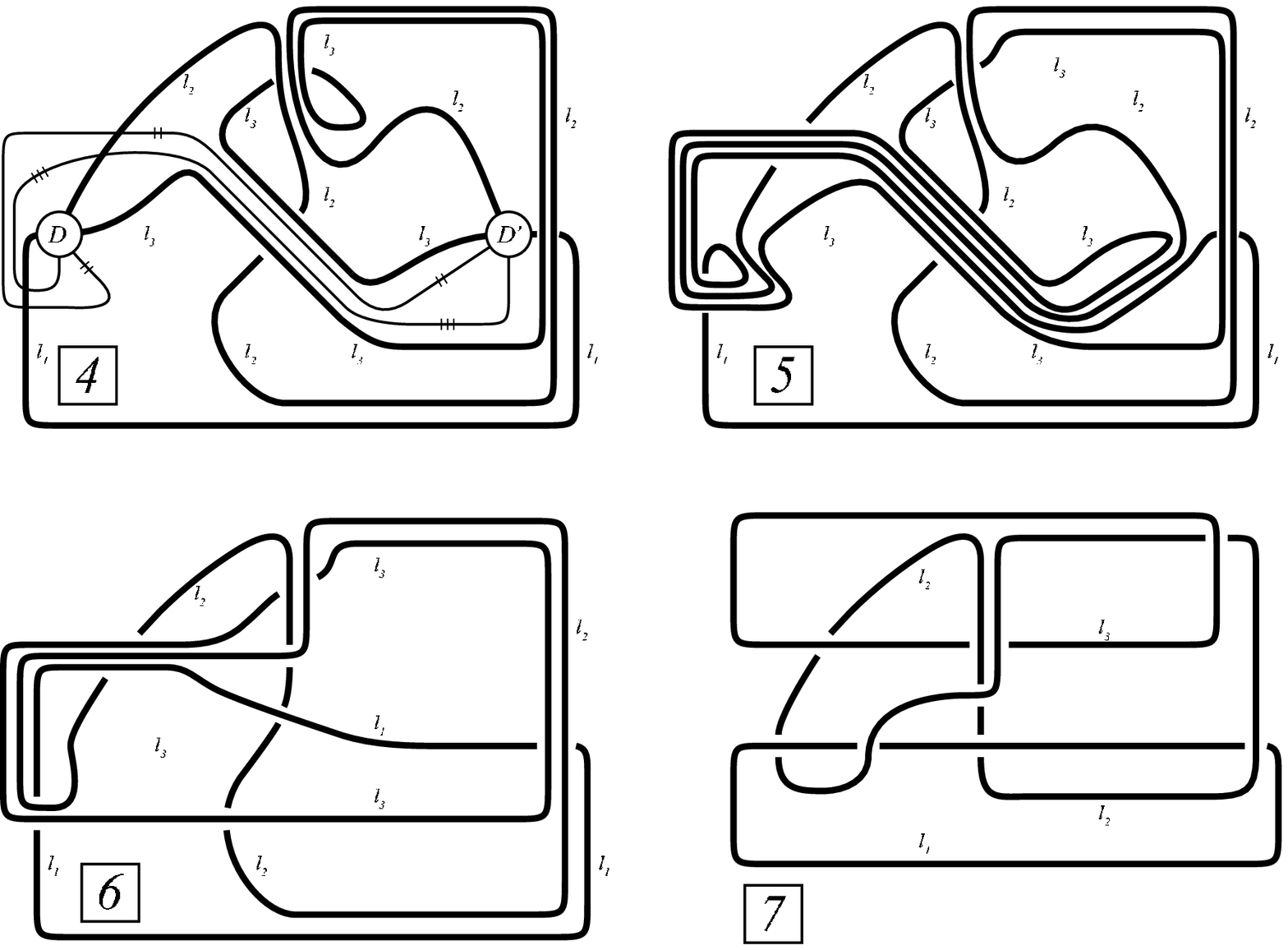}}
\caption{Handle moves, steps 4--7}
\label{handlecancellation2}
\end{center}
\end{figure}

Thus, step~0 of Fig.~\ref{handlecancellation1} shows the handle decomposition of a closed
manifold that we hope is $S^3$.
In the diagrams in Figures~\ref{handlecancellation1} and \ref{handlecancellation2} we perform handle moves
in order to simplify the handle decomposition (see \cite{Gompf-Stipsicz} for basics on handle moves).
Keep in mind that the curves labeled $l_1$, $l_2$ and $l_3$ are not attaching circles,
but merely curves drawn on the surface of the handlebody whose position we keep track of.
In particular, attaching circles may freely be isotoped over these curves and may cross them.
It is easy to see that a crossing by an attaching circle will become an undercrossing if the
corresponding 2-handle cancels a 1-handle that carries one of the longitudes.

{\it Step 0.}
Attaching circles $m_2$ and $m_3$ go across 1-handles $AA'$ and $CC'$ only once, respectively,
so their corresponding 2-handles cancel the 1-handles $AA'$ and $CC'$.
Step~1 shows the handle decomposition after this cancellation.

{\it Step 1.}
Attaching circles II and III, which loop from feet $B'$ and $B$ can be slid over the 1-handle 
$BB'$ and then off feet $B$ and $B'$, respectively.
Moreover, the looping part of attaching circle I, at near right, may be isotoped to foot $D'$ and then across
and off handle $DD'$, after which I is a simple closed curve bounding a disk (on the outside) that may
be pushed away from the diagram.
A 2-handle whose attaching circle bounds a disk disjoint from the rest of the diagram simply encloses
a 3-handle if the manifold is compact, like in our case.
The 2- and 3-handles then cancel.
Step~2 shows the handle decomposition after the isotopies and cancellation.

{\it Step 2.}
We now notice that the vertical portion of attaching circle II can be isotoped outside of the diagram
to the right and ``wrapped'' across $\infty$ to its new position shown in Step~3.
Attaching circle $m_3$ crosses 1-handle $BB'$ only once, causing cancellation of the 2-handle
corresponding to $m_3$ and the 1-handle $BB'$.

{\it Step 3.}
Before we carry out further cancellation, we simplify the picture a bit.
We isotope attaching circle III around $D'$ so it attaches at the bottom.
Notice that this includes a slide of the part of III that runs across the 1-handle $DD'$, thus the place where III
attaches to $D$ moves as well.
Also, we isotope the loop of $l_2$ at the bottom of the diagram  toward the top, and we straighten
out the kink in $l_3$.

{\it Step 4.}
We isotope at top middle to remove the self-crossing of $l_3$.
Attaching circles II and III run parallel, that is, they bound an annulus.
This means a 3-handle is located between them which cancels one of the 2-handles, say III.
After erasing III we note that II cancels the 1-handle $DD'$.

{\it Step 5.}
The rest is isotopy of the link components $l_i$, $i=1,2,3$.
The loop of $l_3$ at center right is isotoped up and to the left, and so is the section of $l_2$ close to it.
The kinks on the left are straightened out, as is the bottom part of $l_3$ and the center of $l_1$.

{\it Step 6.}
The bottom part of $l_3$ is lifted and flipped to the top,
$l_1$ is straightened out and $l_2$ is isotoped a little.

{\it Step 7.}
After isotoping $l_2$ and rotating the diagram by $180^\circ$, one gets the mirror image of Fig.~7 from
Wielenberg's paper~\cite{Wielenberg}.

\section{Converting a side pairing to a handle decomposition in dimension $n$}
\label{convgen}

In this section we generalize the conversion method from~\S\ref{conv3} to any dimension.
Let $X=\en$, $\sn$ or $\hn$, and let $P$ be a finite-volume, finite-sided polyhedron in $X$, as defined,
for example, in \cite{Ratcliffe}.
Assume, furthermore, that every $k$-face of $P$ is diffeomorphic to either $D^k$ or
$D^k-\{\text{finitely many points on $\bd D^k$}\}$.
The former condition is in order to disallow polyhedra such as a lens in $S^3$, whose 1-face is a circle, the
latter to allow hyperbolic polyhedra with ideal vertices.
Let there be given a side-pairing on the sides of $P$, again in the sense of \cite{Ratcliffe}, so that
the space of identified points $M=P/ \sim$ is a complete manifold with a geometry based on $X$.
If $M$ is a noncompact hyperbolic manifold, we will be getting the handle decomposition of $\mbar$, the
compact manifold with boundary whose interior is $M$.
If $M$ is closed, we set $\mbar=M$.

\begin{theorem}
\label{gendecomp}
Let $M$ be a manifold obtained through a side-pairing defined on a polyhedron $P$ in
$X=\en$, $\sn$ or $\hn$.
Suppose that every $k$-face of $P$ is diffeomorphic to either $D^k$ or 
$D^k-\{\text{finitely many points on $\bd D^k$}\}$.
Then the decomposition of $P$ into $k$-faces, $0\le k \le n$ induces a handle decomposition
of the manifold $\mbar$, where every cycle of $k$-faces corresponds to an $(n-k)$-handle.
\end{theorem}

{\it Proof.}
If $P$ is a hyperbolic polyhedron with ideal vertices, the completeness of $M$ implies the existence
of a finite collection of disjoint open horoballs $\{B_s, s\in S\}$ that are centered at ideal vertices of $P$
and are mapped to each other under side-pairings of $P$.
Furthermore, each $B_s$ can be chosen so it intersects only sides of $P$ that are incident with the ideal
vertex where $B_s$ is centered.
Set $U_{-1}=\cup_{s\in S} B_s$ if $P$ is hyperbolic with ideal vertices, otherwise set  $U_{-1}=\emptyset$.

For every $k=0,\dots,n$, we inductively define real numbers $\epsilon_k$ and ``orthogonal neighborhoods''
$NE^k$ of truncated $k$-faces $TE^k$.
Let $E^k_s$, $s\in S_k$ be the collection of  $k$-faces of $P$, $k=0,\dots,n-1$.
(Note that there is only one $n$-face, namely $P$.)
If $P$ has real vertices, there is an $\epsilon$ so that $p: X\to M$ is injective on $\epsilon$-balls
around the real vertices.
Set $\epsilon_0$ to be the smaller of $\epsilon$ and $\frac{1}{3}\min \{ d(E^0_s, E^0_t) | s\ne t, s,t\in S_0\}$, 
otherwise (if all vertices are ideal) set $\epsilon_0=1$.
Let $TE^0_s=E^0_s$, let $NE^0_s$ denote the closed $\epsilon_0$-neighborhhood in $X$ of a 0-face
$E^0_s$, and let $U_0=\cup_{s\in S_0} \intr NE^0_s$.
Clearly $NE^0_s$ and $NE^0_t$ are disjoint when $s\ne t$ and $p$ restricted to any of those
neighborhoods is a diffeormorphism.

Now assume $\epsilon_k$, $U_k$, $NE^k_s$ and $TE^k_s$ have been defined for a $0\le k\le n-2$ and
every $k$-face $E^k_s$, $s\in S_k$ of $P$ and that the restriction of $p$ to every $NE^k_s$ is a
diffeomorphism.
Let $TE^{k+1}_s=E^{k+1}_s - \cup_{-1 \le i \le k}  U_i $, $s\in S_{k+1}$.
Because $M$ is a manifold, $p$ is injective on the interior of every face.
Due to compactness of every $TE^{k+1}_s\subset \intr E^{k+1}_s$, we can find an $\epsilon$ so that 
$p$ is a diffeomorphism on an $\epsilon$-neighborhood of $TE^{k+1}_s$ for every $s\in S_{k+1}$.
Set $\epsilon_{k+1}$ to be the smallest of $\epsilon$, $\frac{1}{2}\epsilon_k$ and
$\frac{1}{3}\min \{ d(TE^{k+1}_s, TE^{k+1}_t) | s\ne t, s,t\in S_{k+1}\}$,  and let $NE^{k+1}_s$ be the closed
$\epsilon_{k+1}$-neighborhood of $TE^{k+1}_s$ in $X$ with $\cup_{i=-1}^k  U_i $
excluded.
Let $U_{k+1}=\cup_{s\in S_{k+1}} \intr NE^{k+1}_s$.
If $k=n-1$, let  $TE^n=P-\cup_{i=-1}^{n-1} U_i $ and $NE^n=TE^n$.

From the assumption that every $E^k_s$ is diffeomorphic to $D^k$ or $D^k-\{\text{finite set}\}$, it
follows that $TE^k_s$ is diffeomorphic to $D^k$ for every $s\in S_k$.
The set $NE^k_s$ is then diffeomorphic to $D^k\times D^{n-k}$, where $TE^k_s=D^k\times 0$.
Note that $x\times D^{n-k}$ is essentially in the orthogonal direction to $E^k_s$, except close to
$\bd TE^k_s$, where some bending has to occur to accomodate $NE^i_t$, where $E^i_t$ is a face of $E^k_s$.
Furthermore, for every face $E^i_t$ of $E^k_s$, $i\le k$, note that $NE^i_t$ intersects
$NE^k_s=D^k\times D^{n-k}$ only along $\bd D^k\times D^{n-k}$.
If $P$ has ideal vertices, $\bd\mbar$ is assembled from links of ideal vertices $\bd B_s \cap P$.
Clearly $\bd B_s \cap P\subset \bd D^k\times D^{n-k}$.
Thus, any element of $p(\intr D^k \times\bd D^{n-k})$ is in $\intr \mbar$ and therefore must be in some other
$p(NE^j_u)$.
Our observation above then shows that $j>k$.

Let us treat $p(NE^n)$ as a $0$-handle in $\mbar$.
Consider an $NE^k_s$, $s\in S_k$, $k\le n$.
By our construction, $p$ restricts to a bijection on $NE^k_s$, hence $p(NE^k_s)$
is an $n$-ball inside of $\mbar$.
This gives a decomposition of $\mbar$ into a collection of $n$-balls with disjoint interiors.
If $E^k_s$ and $E^k_t$ are in the same cycle, then $p(NE^k_s)=p(NE^k_t)$, hence
every $n$-ball corresponds to a cycle of $k$-faces for some $k\le n$.

Define $M_k=\cup_{n-k \le i \le n} \cup_{s\in S_i}  p(NE^i_s)$, meant
to be the union of $i$-handles, $0\le i \le k$.
As  above, $NE^k_s=D^k\times D^{n-k}$ and $p(\intr D^k\times D^{n-k})$ is contained
in $\cup_{k < i} \cup_{t\in S_i} p(NE^i_t)=M_{n-k-1}$.
Since $M_{n-k-1}$ is closed, $p(D^k\times \bd D^{n-k})\subset M_{n-k-1}$.
Therefore, $p(NE^k_s)$ attaches as an $(n-k)$-handle to $M_{n-k-1}$, giving us a handle
decomposition of $\mbar$.
$\qed$

In the above handle decomposition of $\mbar$, we note that the attaching sphere $0\times \bd D^{n-k}$
of the $n-k$-handle $NE^k_s$ is the boundary of a neighborhood in $X$ of a point
$x\in E^k_s$.
Naturally, this being a neighborhood in $X$ means that a part of it is outside of $P$ and it intersects
several translates $gP$ of $P$.
But then $g^{-1}((0\times\bd D^{n-k}) \cap gP)$ is visible in $P$.

\section{Drawing handle decomposition diagrams in dimension 4}
\label{diagram4}

In this section we apply the conversion method described in the previous section to dimension~4.
Notation is like in the previous section and, as an illustrative example, we use as $P$ the 4-cube
whose sides are paired by translations, yielding $M=$4-torus.

We want to draw in $\bd D^4=S^3=\riii \cup \infty$ attaching spheres of the $k$-handle
$D^{4-k}\times D^k$.
The 0-handle is $NE^4$, $P$ without neighborhoods of all the $k$-faces.
Clearly $\bd NE^4 =S^3=\bd P$, realized by a diffeomorphism $h:\bd NE^4\to \bd P$, 
a restriction of a diffeomorphism $h:NE^4\to P$, which may be imagined as a radial projection
from a point in the interior of $P$.
Under~$h$, $(\bd NE^{4-k}_s)\cap P$ is sent to $TE^{4-k}_s\times B^{k-1}$, a $TE^{4-k}_s$
``thickened-up''   in $\bd P$.
Note that the thickening of $TE^3_s$ in $\bd P$ is still $TE^3_s$.
Now, a piece of the attaching sphere $P\cap (0\times \bd D^k)$ is sent under $h$ to 
$x \times B^{k-1}$, where $x\in TE^{4-k}_s$.

Let the subdivision of $\bd P$ into $k$-faces ($k\le 3$) be drawn in $\riii\cup\infty$.
As an example, take the standard ``cube-within-a-cube'' picture of the boundary of the 4-cube.
A piece of the attaching sphere for a 1-handle $p(NE^3_s)$ is $P\cap (0\times \bd D^1)$ which is
sent to a point in $TE^3_s$, chosen, for example, in its interior.
The two points in the attaching sphere of a 1-handle are in paired 3-faces.
The attaching region $D^3\times \bd D^1$ is the union of paired truncated 3-faces $TE^3_s$ and
$TE^3_t$: schematically, we draw 3-balls inside $TE^3_s$ and $TE^3_t$.

\begin{figure}
\resizebox{\textwidth}{!}{\includegraphics{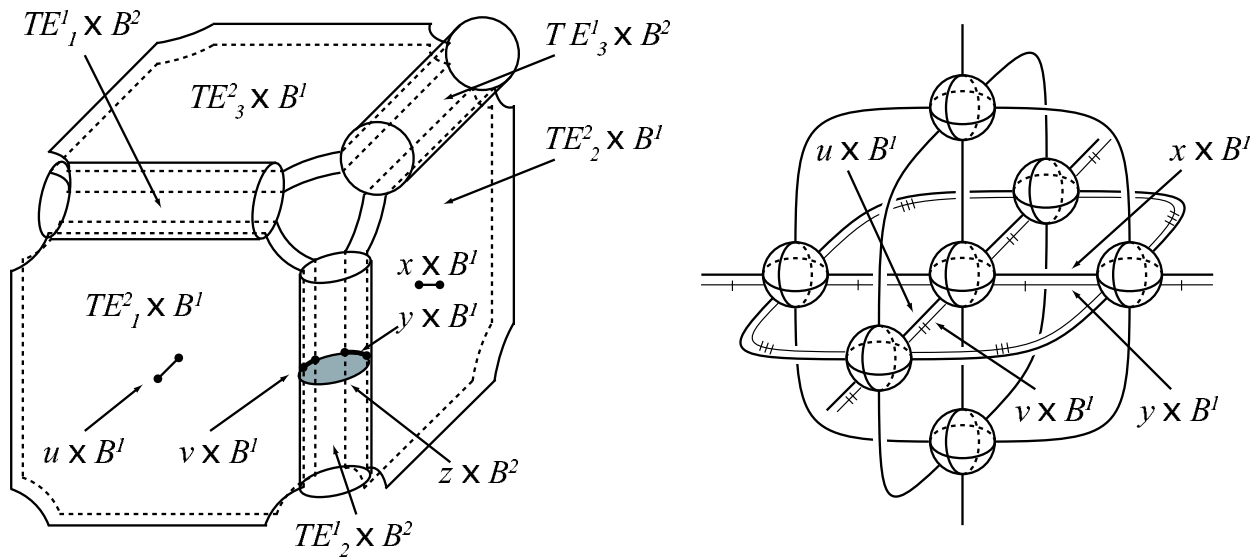}}
\caption{Arriving at a handle decomposition for the 4-torus}
\label{dim4conv}
\end{figure}

The piece inside of $P$ of an attaching circle for a 2-handle $p(NE^2_s)$ is $P\cap (0\times \bd D^2)$,
an arc that crosses the 2-face corresponding to the 2-handle and joins the 3-faces whose intersection
is the 2-face.
Under $h$, this arc maps to the segment $x\times B^1$, $x\in TE^2_s$, visible on the left of Fig.~\ref{dim4conv}.

The attaching sphere for a 3-handle is a 2-sphere, whose intersection $P\cap (0\times \bd D^3)$
with $P$ is a 2-ball.
Under $h$, the 2-ball maps to the 2-ball $z\times B^2$, $z\in TE^1_s$, shown in Fig.~\ref{dim4conv}.

When $M$ is closed, it is only important how the 1- and 2-handles attach (see \cite{Gompf-Stipsicz}),
so the attaching spheres of 3- and 4-handles do not matter.
However, it is useful to have in one's mind pieces of the attaching spheres of the 3-handles,
as they help with the framing of the attaching map of the 2-handles.

In order to specify, up to isotopy, the attaching map $\phi:D^2\times \bd D^2\to \bd Y$ of a 2-handle
$D^2\times D^2$, it is enough to specify the images of two parallel circles $\phi(x\times \bd D^2)$
and $\phi(y\times \bd D^2)$.
As we can see on the left side of Fig.~\ref{dim4conv}, if $E^1_t$ is incident with $E^2_s$, then the intersection
of $z\times B^2$, $z\in TE^1_t$, with $TE^2_s\times B^1$ is an arc $y\times B^1$, $y\in TE^2_s$.
We can then choose $y\times \bd D^2$ to be the circle parallel to $x\times \bd D^2$, chosen before.
(In Fig.~\ref{dim4conv}, $u\times B^1$ and $v\times B^1$ are pieces of another such pair of parallel circles.)

Schematically, the portion $z\times B^2$ of the attaching sphere of the 3-handle $p(NE^1_t)$ is
represented as a triangle transverse to $E^1_t$, bounded by the portions $x_i \times B^1$
of attaching circles of the three 2-handles that correspond to the three 2-faces, the
pairwise intersections of the three 3-faces whose intersection is the 1-face $E^1_t$.
One can speak of a ``cycle'' of triangles, the collection of triangles corresponding
to 1-faces that are all in one cycle.
Clearly, a cycle of triangles represents all the pieces of the attaching sphere of the 3-handle
corresponding to the cycle of 1-faces.

Thus, pieces of a parallel circle may be chosen to lie in one of the triangles.
Since a 2-face is incindent with several 1-faces, pieces of the attaching circle of a 2-handle
will be in the boundary of several triangles.
We choose one to contain a piece of the parallel circle; once this is done, the remaining pieces of the
parallel circle must be chosen in triangles that are in the same cycle as the one we have chosen.

We summarize how to get a picture of a handle decomposition of a 4-manifold that is the result
of pairing sides of a polyhedron $P$.

\begin{itemize}
\item[---]
Draw in $\riii$ the decomposition of $\bd P=\riii\cup \infty$ into $k$-faces .
Inside every 3-face, draw a 3-ball.
Feet of a 1-handle are the two balls inside paired 3-faces.
\item[---]
We do not assume that the feet of 1-handles are identifed by a reflection in the
bisector of the centers, as is common.
Rather, the identifying map is determined by the map that pairs the corresponding 3-faces.
\item[---]
If two 3-faces are adjacent along a 2-face $E^2$, draw an arc between the balls inside
the 3-faces that crosses $E^2$ exactly once.
The arcs that cross 2-faces that are in the same cycle comprise the attaching circle for a 2-handle.
\item[---]
Whenever three 3-faces intersect in a 1-face we see a ``triangle'' whose ``vertices'' and edges are
the already drawn 3-balls and arcs, respectively.
We fill in this triangle (usually only mentally) with a surface that is transverse to the 1-face.
Parallel attaching circles can be chosen to lie in these surfaces.
\item[---]
Once we choose a triangle to contain a piece of the attaching circle, the remaining pieces
must be chosen in triangles that are in the same cycle of triangles.
\end{itemize}

The procedure above yields the familiar handle-decomposition diagram for $T^4$ from
the right side of Fig.~\ref{dim4conv}.
(see also \cite{Gompf-Stipsicz}, Fig.~4.42).
Parallel attaching circles for three 2-handles are the arcs marked I, II and III.

\section{A hyperbolic manifold as a complement of 5 tori in the standard differentiable $S^4$}
\label{ident4}

In \cite{Ivansic3}, the author showed that the double cover of $M_{1011}$, example no. 1011 from
Ratcliffe and Tschantz's  \cite{Ratcliffe-Tschantz} collection of noncompact finite-volume
hyperbolic 4-manifolds, is a complement of 5 tori in the topological 4-sphere.
The proof used Freedman's theory, which only provides a homeomorphism to the 4-sphere.
In this section we prove that the 4-sphere $N$ from \cite{Ivansic3} is, in fact, diffeomorphic to the
standard differentiable 4-sphere.
We use the method of this paper to obtain a handle decomposition of the manifold $N$ and then
handle moves to simplify the decomposition down to the decomposition of the standard differentiable
4-sphere.

The 24-sided polyhedron $Q$ that gives rise to the Ratcliffe and Tschantz's manifolds is described
in their paper \cite{Ratcliffe-Tschantz} and in \cite{Ivansic3} and \cite{Ivansic4}, where more details
on its combinatorial structure can be found.
Here we just recall that its sides (in the ball model of $\hiv$) are spheres of radius 1 centered at
points whose two coordinates are $\pm 1$ and the other two are zeroes.
We label the spheres and the sides by $S_{****}$, like in \cite{Ivansic4}. 
For example, $S_{0+0-}$ is the sphere centered at $(0,1,0,-1)$.

Each octahedral 3-face of $Q$ has eight 2-faces, so drawing a decomposition of $\bd Q$ would be quite
involved.
We will therefore jump to the handle-decomposition picture right away, by finding
attaching spheres of the 1-, 2- and 3-handles on $\bd Q$, and projecting them to $S^3$ radially from
the origin of $B^4$.
$S^3$ is then sent to $\riii\cup \infty$ via the M\"obius transformation $g:(\riv\cup\infty)\to(\riv\cup\infty)$
that provides the standard isometry between the ball and upper-half-space models of hyperbolic space.
This map is the composite of the reflection in the sphere with center $(0,0,0,1)$ of radius $\sqrt 2$,
followed by a reflection in the hyperplane $x_4=0$.
Its restriction to $S^3$ is given by  $x\mapsto e_4 + \frac{2}{|x-e_4|^2}(x-e_4)$.
(This is actually the formula for just the first reflection, since the reflection in $x_4=0$ has no effect on
$\riii\cup\infty$, the image of $S^3$.)
Note that $g$ leaves $S^2\subset \riii\times 0$ fixed.

As attaching spheres of 1-handles we choose points on the sides of $Q$ closest to the origin.
``Shortest'' arcs connecting those points along the sides are chosen to be the pieces of attaching circles
of 2-handles.
Pieces of attaching spheres of 3-handles are the piecewise-spherical ``triangles'' bounded
by the arcs, stretched across the sides.

More precisely, let $r$ be the position-vector of a sphere $S$ that determines a side of $Q$.
The intersection of $S$  and the line spanned by $r$ is a point $c$ in $S$.
Let $c'$ be the intersection of $S'$ and the line spanned by $r'$, where $S'$ is the pair of $S$
under the side-pairing on $Q$.
Then we choose $c$ and $c'$ to be the points of the attaching sphere of the 1-handle corresponding
to the paired sides $S$ and $S'$.
If $S_1$ and $S_2$ are intersecting sides, let $e$ be the arc that is the intersection of $\bd Q$
and the (linear) angle spanned by position vectors $r_1$ and $r_2$.
This arc is a portion of the attaching circle of the 2-handle corresponding to the 2-face $S_1\cap S_2$.
Finally, if a side $S_3$ intersects $S_1$ and $S_2$, consider the intersection $f$ of $\bd Q$
with the ``positive cone'' spanned by $r_1$, $r_2$ and $r_3$.
This ``triangle'' is a portion of the attaching sphere of the 3-handle corresponding to the
1-face $S_1\cap S_2 \cap S_3$.
It is not clear that the overall arrangement of the spheres is such that $c$, $c'$, $e$ and $f$ are
actually on $\bd Q$.
(For example,  $c$ may be inside some other sphere $S_0$, which would put it outside of $Q$.)
Next, we justify that these choices are indeed on $\bd Q$.

Consider the sides $S_{++00}$, $S_{+0+0}$ and $S_{0++0}$.
The three sides intersect pairwise and the intersection of all three of them is a 1-face $E^1_s$.
Let $L$ be the (linear) hyperplane spanned by $r_{++00}$, $r_{+0+0}$ and $r_{0++0}$, this is the
hyperplane $x_4=0$.
It is clear that the only spheres that intersect $L$ in more than one point are those with a 0 in the
fourth position of its label.
All other spheres intersect $L$ in exactly one point, which is one of $\pm e_i$, $i=1,..,4$, that is, 
an ideal vertex of $Q$.
Now, the intersection of  the 12 sides $S_{***0}$ with $L$ is a 3-dimensional version of $Q$
which was described and pictured in \cite{Ratcliffe-Tschantz}, Figure 5.
From this picture we see that attaching spheres chosen in the way described above are on $\bd Q$.

\begin{figure}
\resizebox{\textwidth}{!}{\includegraphics{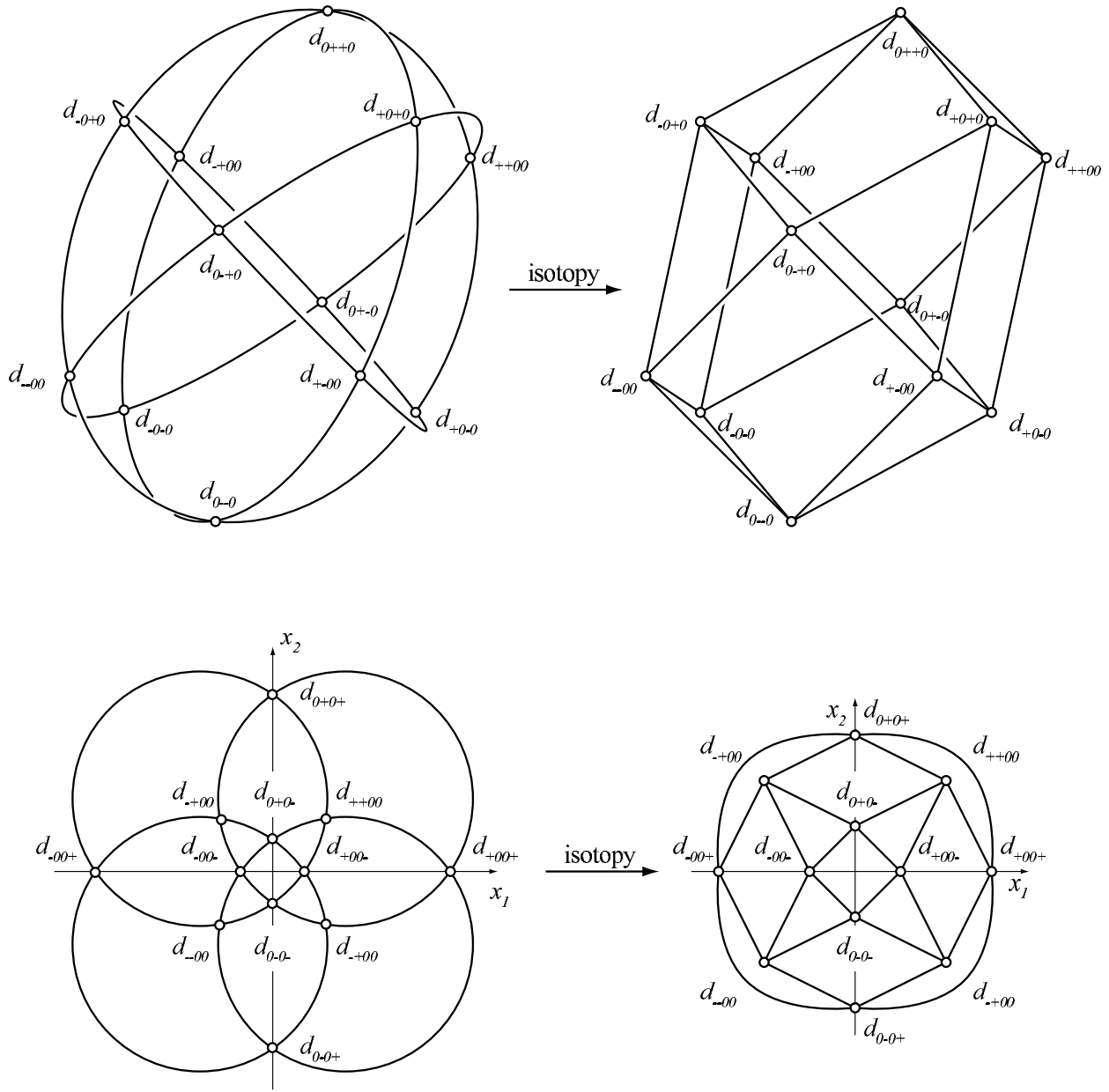}}
\caption{Finding attaching spheres for $M_{1011}$}
\label{polytohandle}
\end{figure}

A general 1-face is an intersection of 3 sides if their labels pairwise share exactly one position
with the same symbol.
This position could be different for each pair of sides, like in the example above, or it could be the
same for all three pairs, like for the sides $S_{++00}$, $S_{+0+0}$ and $S_{+00+}$.
It is clear that any 1-face can be moved by a linear isometry of $Q$ to one of these two
prototypical 1-faces (permute the coordinates and reflect in coordinate hyperplanes).
Furthermore, there is a linear isometry of $Q$ that sends 
$S_{++00}$, $S_{+0+0}$ and $S_{0++0}$ to $S_{++00}$, $S_{+00+}$ and $S_{+0+0}$, respectively;
its matrix is
\begin{displaymath}
\frac{1}{2}
\left[
\begin{array}{rrrr}
1 & 1 & 1 & -1\\
1 & 1 & -1 & 1\\
-1 & 1 & 1 &1\\
1 & -1 & 1 & 1
\end{array}
\right].
\end{displaymath}
This shows that the situation illustrated by the sides $S_{++00}$, $S_{+0+0}$ and $S_{0++0}$
is generic, so all choices for attaching spheres done in the way described above are valid.

We now have to see where the attaching spheres are sent under the composite $gp $,
where $p:\bd Q\to S^3$ is the radial projection.
The intersection of each position vector $r_{****}$ with $S^3$ is $\frac{1}{\sqrt 2}r_{****}$.
The points of form $\frac{1}{\sqrt 2}r_{***0}$ are on $S^2\subset S^3$, which is fixed by $g$.
Furthermore, an easy computation shows that $g(\frac{1}{\sqrt 2}r_{***+})=(\sqrt 2+1)(*,*,*,0)$ and
$g(\frac{1}{\sqrt 2}r_{***-})=(\sqrt 2-1)(*,*,*,0)$.

As above, let $c_i$ be the point of intersection of a sphere (side) $S_i$ with the line spanned by the position
vector $r_i$ of the center of $S_i$, $i=1,2,3$.
If sides $S_1$ and $S_2$ intersect, consider the intersection $C$ of $S^3$ and the linear plane spanned
by $r_1$ and $r_2$ (part of $C$ is the radial projection of a piece of the attaching circle).
This is a circle, so $g(C)$ is a circle, since $g$ is a M\"obius transformation.
Since $C$ also contains $-r_1$ and $-r_2$, the circle $g(C)$ will contain the four points $gp(\pm c_1)$
and $gp(\pm c_2)$.
Once we have the four points drawn, the circle $g(C)\subset\riii$  will be easy identify.
The arc of the circle between $gp(c_1)$ and $gp(c_2)$ is a part of the attaching circle for the 2-handle
corresponding to the 2-face $S_1\cap S_2$.
Now, if sides $S_1$, $S_2$ and $S_3 $ all intersect in a 1-face $E^1$, part of the attaching sphere
corresponding to $E^1$ is the ``triangle''  $f$ that is the intersection of the positive cone generated
by $r_1$, $r_2$ and $r_3$ with $\bd Q$.
Radial projection to $S^3$ followed by $g$ maps the triangle to a spherical triangle bounded by arcs of circles
$g(C)$ that were just described.

The top left  of Fig.~\ref{polytohandle} shows the points $d_{****}=gp(c_{****})$ and the circles $g(C)$
for the sides $S_{**00}$, $S_{*0*0}$ and $S_{0**0}$.
The bottom left of Fig.~\ref{polytohandle} does the same for the sides  $S_{**00}$, $S_{0*0*}$ and $S_{*00*}$.
The complete picture for $Q$ is obtained by rotating the bottom left by $\pi/2$ around the $x_1$-axis,
then around the $x_3$-axis, and taking the union of the resulting three figures with the
top left of Fig.~\ref{polytohandle}.
To make drawing of pictures easier, we isotope the positions of $gp(c)$'s a little and replace curved
arcs $g(C)$ mostly by straight lines, as seen on the right side of Fig.~\ref{polytohandle}.

We note that pieces of the attaching circles all lie in one of the coordinate planes or on~$S^2$.
(In the straight-edge version of the diagram we imagine this $S^2$ as the surface consisting of 
6 rectangles and 8 triangles, spanned by the points $g(c)$).
We assume, as discussed in \S\ref{diagram4}, that pieces parallel circles always lie in the triangles
in the diagram.
Note that each triangle can be taken to lie in one of the coordinate planes or on $S^2$.

\begin{figure}
\resizebox{\textwidth}{!}{\includegraphics{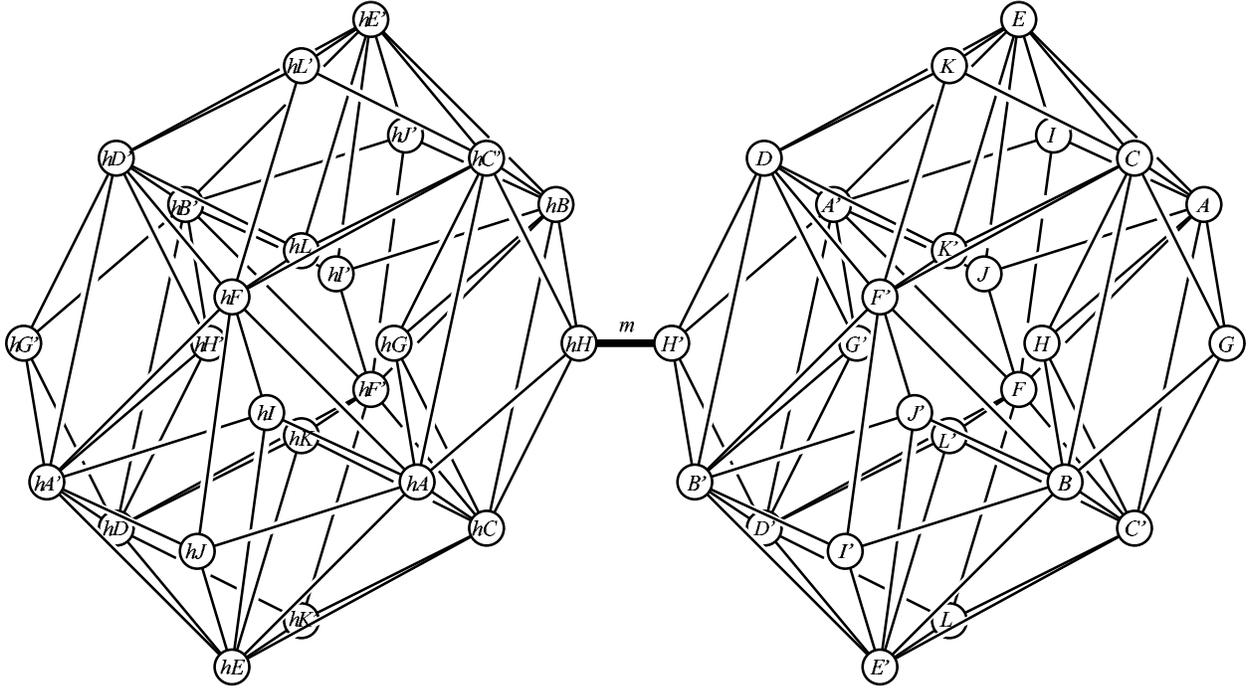}}
\caption{Handle decomposition of $\tilde M_{1011}$}
\label{handle1011}
\end{figure}

\begin{figure}
\begin{center}
\resizebox{1.5in}{!}{\includegraphics{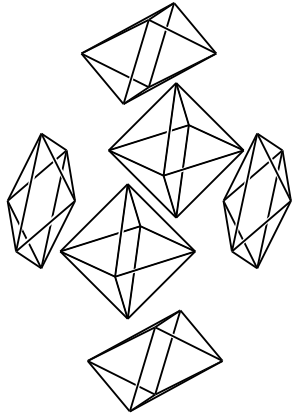}}
\caption{Octahedra appearing in Fig.~12, separated}
\label{octahedra}
\end{center}
\end{figure}

The handle decomposition of $\mbar_{1011}$ is the right half of Fig.~\ref{handle1011}, where
the outside- and inside-most attaching circles are not shown to maintain clarity of picture.
If every triangle in the picture is filled in, we note that it will consist of eight ``octahehedra", each of
which corresponds to the link of the ideal vertices $v_{*000}$, $v_{0*00}$, $v_{00*0}$ and
$v_{000*}$, which is a cube (three pairs of opposing sides/feet of 1-handles).
To better see the octahedra, we have separated them on Fig.~\ref{octahedra}.
Note that six are shown; the two missing ones, corresponding to $v_{000+}$ and $v_{000-}$, are the same
ones that are missing from Fig.~\ref{handle1011}.
As a matter of fact,  the space between the described octahedra forms sixteen more octahedra,
corresponding to the ideal vertices of form $v_{****}$ (see one in Fig.~\ref{meridian5}).

A side-pairing $f:S\to S'$ of any of Ratcliffe and Tschantz's examples is always of the form $ru$,
where $u$ is a composite of reflections in the coordinate planes, and $r$ is a reflection in $S'$.
The restriction of $f$ to $\bd Q$ is all that matters to us, so $u$ explains how feet of 1-handles are
identified.
Note that conjugating by $g$ the reflections in hyperplanes $x_1=0$, $x_2=0$ and $x_3=0$ gives the 
same reflections, while conjugating the reflection in $x_4=0$ gives the reflection in the unit sphere.
Using the convention from \cite{Ivansic4} we name the side-pairings for $M_{1011}$ by letters
$a,b,\dots,k, l$ as follows.
(The composite of reflections that pair the sides is under the arrow.)
\begin{displaymath}
\begin{array}{llll}
S_{++00} \arrtop{a}{-+++} S_{-+00}\hskip10pt
&
S_{+-00} \arrtop{b}{-+++} S_{--00}\hskip10pt
&
S_{+0+0} \arrtop{c}{++-+} S_{+0-0}\hskip10pt
&
S_{-0+0} \arrtop{d}{++-+} S_{-0-0}
\\
S_{0++0} \arrtop{e}{----} S_{0--0}\hskip10pt
&
S_{0+-0} \arrtop{f}{----} S_{0-+0}\hskip10pt
&
S_{+00+} \arrtop{g}{----} S_{-00-}\hskip10pt
&
S_{+00-} \arrtop{h}{----} S_{-00+}
\\
S_{0+0+} \arrtop{i}{+-++} S_{0-0+}\hskip10pt
&
S_{0+0-} \arrtop{j}{+-++} S_{0-0-} \hskip10pt
&
S_{00++} \arrtop{k}{+++-} S_{00+-}\hskip10pt
&
S_{00-+} \arrtop{l}{+++-} S_{00--}.
\end{array}
\end{displaymath}
Furthermore, for simplicity of notation, if a letter $s$ pairs two sides,  we relabel the originating side
with $S$ and $s(S)$ by $S'$.
Thus, $d=$reflection in plane $x_3=0$ sends side $D=S_{-0+0}$ to side $D'=S_{-0-0}$.

Let $G_{1011}\subset \Isom(\hiv)$ be the fundamental group of $M_{1011}$ and let $H_{1011}$
be the subgroup of orientation-preserving isometries in $G_{1011}$.

Of course, $G_{1011}$ is generated by $a,b,\dots,l$.
We are really interested in the orientable double cover $\tilde M_{1011}$ of $M_{1011}$,
whose fundamental polyhedron consists of two copies of $Q$ with suitably paired sides.
It is easy to see (and is explained in \S3 of \cite{Ivansic3}) that the fundamental polyhedron
for $H_{1011}$ is $Q\cup hQ$, where $h$ is one of the above listed generators of $G$, that is also,
being orientation reversing, a  coset representative for the nontrivial right coset of
$H_{1011}$ in $G_{1011}$.
The discussion in \cite{Ivansic3} also shows that sides of $Q\cup hQ$ are paired according to
the following rule.
Let $S$, $S'$ be sides of $P$ paired by the transformation $s\in G$.
If $s$ is orientation-reversing, then side $S$ is paired to $hS'$ via $hs$ and side $hS$ is paired to
$S'$ via $sh^{-1}$.
If $s$ is orientation-preserving, then $S$ is paired to $S'$ via $s$ and $hS$ is paired to $hS'$
via $hsh^{-1}$.

We may view the handle decomposition of $\tilde M_{1011}$ as having two 0-handles (corresponding
to $Q$ and $hQ$) and a 1-handle joining them (coming from the paired sides $H'$ and $hH$).
Since $Q$ and $hQ$ lie on opposite sides of the hyperplane $H'$, it is clear that the handle
decomposition of $\tilde M_{1011}$ can be drawn by drawing two handle-decompositions of
$\mbar_{1011}$ side-by-side while identifying $H'$ and $hH$.
The same effect is achieved by drawing them side-by-side and introducing a 2-handle canceling
the 1-handle coming from pairing $H'$ to $hH$.

The  handle decomposition for $\tilde M_{1011}$ is the entire diagram in Fig.~\ref{handle1011}.
The part coming from $Q$ we take to be centered at 0,  the part coming from $hQ$ we center
at $(-6,0,0)$, the two portions being symmetric in the plane $x_1=-3$.
To get proper labeling on the feet of 1-handles of the $hQ$-part, we recall that they are the result
of applying $h=ru_{----}$ to $Q$, thus we need to reflect the picture on the right in planes $x_1=0$,
$x_2=0$, $x_3=0$ and the unit sphere centered at 0, and then apply the reflection $r$ in the plane
$x_1=-3$.

Putting together all the facts from above, the feet of the 1-handles in the decomposition in Fig.~\ref{handle1011}
have the following identification pattern: 
\begin{gather*}
A, B, C, D, I, J, K, L, hA, hB, hC, hD, hI, hJ, hK, hL \mapsto\\
A', B', C', D', I', J', K', L', hA', hB', hC', hD', hI', hJ', hK', hL'\\
\text{via reflection in the bisector of the feet}\\
E, F, G, H, E', F', G', H' \mapsto hE', hF', hG', hH', hE, hF, hG, hH\\ 
\text{via reflection in $x_1=-3$.}
\end{gather*}
The manifold $\tilde M_{1011}$ has 5 three-torus boundary components, each of those an
$S^1$-fiber bundle over $T^2$.
Closing off the boundary components involves filling in each fiber with a disc, resulting in a closed
manifold $N$.
Equivalently, this can be achieved by attaching a $T^2\times D^2$ to each component of $\bd \tilde M_{1011}$.
A handle decomposition for $T^2\times D^2$ derived from the simplest handle decomposition for $T^2$
has one 0-handle, two 1-handles and one 2-handle.
Attaching it to $\tilde M_{1011}$ results in adding one 2-handle, two 3-handles and one 4-handle
to the decomposition, since the handles in the decomposition of $T^2\times D^2$ must be viewed
in an upside-down way (see \cite{Gompf-Stipsicz}), as it is attached to $\tilde M_{1011}$.
The attaching circle of the 2-handle is any fiber in the bundle.
As selected and illustrated in \cite{Ivansic3}, in four of the boundary components, the fibers are
represented by a straight-line segment joining opposing sides of the cube that is the vertex link,
therefore, the attaching circle of the 2-handle is a line-segment joining two opposed feet of 1-handles
in the octahedron that corresponds to the vertex link.
Since the parallel circle is another fiber, we may assume it is simply a parallel line-segment.

We now simplify the handle decomposition of $N$ using handle moves.
We repeatedly make use of the following proposition:

\begin{proposition}
\label{cancel3handle}
(\cite{Gompf-Stipsicz}, modified Proposition 5.1.9)
If the handle decomposition of a closed manifold contains an attaching circle of a 2-handle that
can be isotoped so that it bounds a disc disjoint from the rest of the diagram, and the disc contains its
parallel circle, then the 2-handle cancels a 3-handle from the decomposition (and we may erase
the 2-handle from the diagram).
\end{proposition}

In order to better see what goes on in the complicated diagram we consider sections
with  the coordinate planes $x_1x_2$, $x_1x_3$, $x_2x_3$ and its parallel plane $x_1=-6$
(``$x_2x_3$-planes''), and the two spheres $S^2$, displayed in this order in
Figures~\ref{step0}--\ref{step6}.
Since pieces of the parallel circles are on one of those surfaces, if isotopy of the diagram
stays parallel to the surface, those pieces remain in the surface, so are easy to track.
On occasion, pieces of parallel circles are not on the default surfaces --- it is either obvious where they are,
or it is noted.

\begin{figure}
\begin{center}
\includegraphics{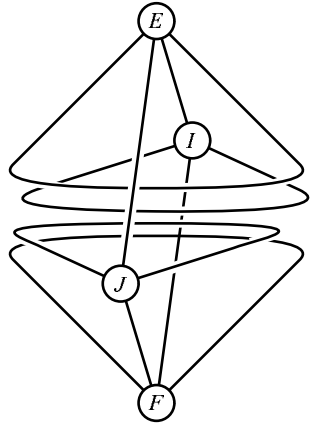}
\caption{Octahedron corresponding to $v_{0+00}$ after $AA'$ is canceled by a 2-handle.}
\label{octahedroncancellation}
\end{center}
\end{figure}

Handle moves are tracked in Figures~\ref{step0}--\ref{step6} by drawing what happens
in the sections with the mentioned planes.
The topmost box of the explanation describes which 1- and 2-handles have canceled.
The middle four boxes describe subsequent isotopies, each of the four
pertaining to the part of the diagram in the same relative position as the box.
The bottom box describes cancellation of 2- and 3-handles owing to
Proposition~\ref{cancel3handle}.
Note that 1-handles are designated by labels on the corresponding paired sides.

One could wonder if isotopy or handle cancellation in one of the planes interferes with the situation
in the other.
A little 3-dimensional insight helps one see that there is no problem.
A picture such as Fig.~\ref{octahedroncancellation}, which is typical, may help the reader
see what happens after a 1-handle cancels with a 2-handle.
This picture shows the octahedron corresponding to $v_{0+00}$ after 1-handle $AA'$ was
canceled by a 2-handle.

The initial handle decomposition also includes thirty-four 3-handles and five 4-handles.
Twenty-four of the 3-handles come from cycles of 1-faces, the remaining ten 3-handles and five
4-handles come from handle decompositions of attached $T^2\times D^2$'s.
As explained in \cite{Gompf-Stipsicz}, \S4.4, because $N$ is a closed manifold, 3- and 4-handles
can attach essentially in only one way, so there is no need to keep track of them.

\begin{figure}
\begin{center}
\includegraphics{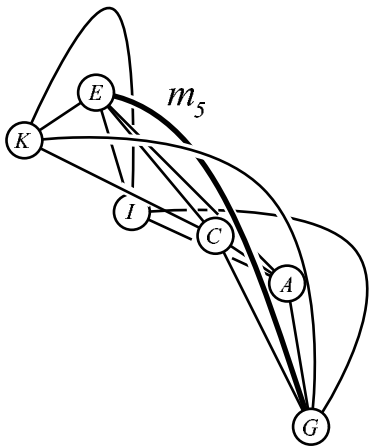}
\caption{Position of $m_5$ in octahedron corresponding to $v_{++++}$}
\label{meridian5}
\end{center}
\end{figure}

After Fig.~\ref{step6}, the 3-dimensional diagram is simple enough that we can draw it on one picture.
In steps thus far, we have not pictured the additional 2-handle coming from closing off the fifth boundary
component.
The choice $e^{-1}g$, made in \cite{Ivansic3}, is represented by a union of line-segments:
one joining the opposite sides $S_{+00+}$ and $S_{0++0}$ of the cube that is the
vertex link at $v_{++++}$, and one joining the opposite sides $S_{0--0}$ and $S_{-00-}$
of the cube that is the vertex link at $v_{----}$.
This corresponds to two segments joining $E$ and $G$ and $hE'$ and $hG'$ in our diagram.
Fig.~\ref{meridian5} shows the first one in the ``octahedron" corresponding to $v_{++++}$: we can
see that all the moves done so far on the diagram do not affect it (in particular, it lies outside of the
two $S^2$'s), so it can be drawn in the same position in the overall picture.

After the final 2-handles and 3-handles cancel in step~12 of Fig.~\ref{step7}, we arrive at an
empty diagram (one 0-handle and some 3- and 4-handles).
Since this is the handle decomposition of the standard differentiable 4-sphere, we conclude
that $N$ is diffeomorphic to it.

Incindentally, by keeping track of the 3- and 4-handles throughout the computation we see that the
final handle decomposition has four 3-handles and five 4-handles.
However, since the boundary of their union is $S^3$, the union must be $D^4$ (otherwise the
boundary would be a connected sum of $S^1\times S^2$'s).

\pagebreak

\begin{landscape}

\begin{figure}
\resizebox{8in}{!}{\includegraphics{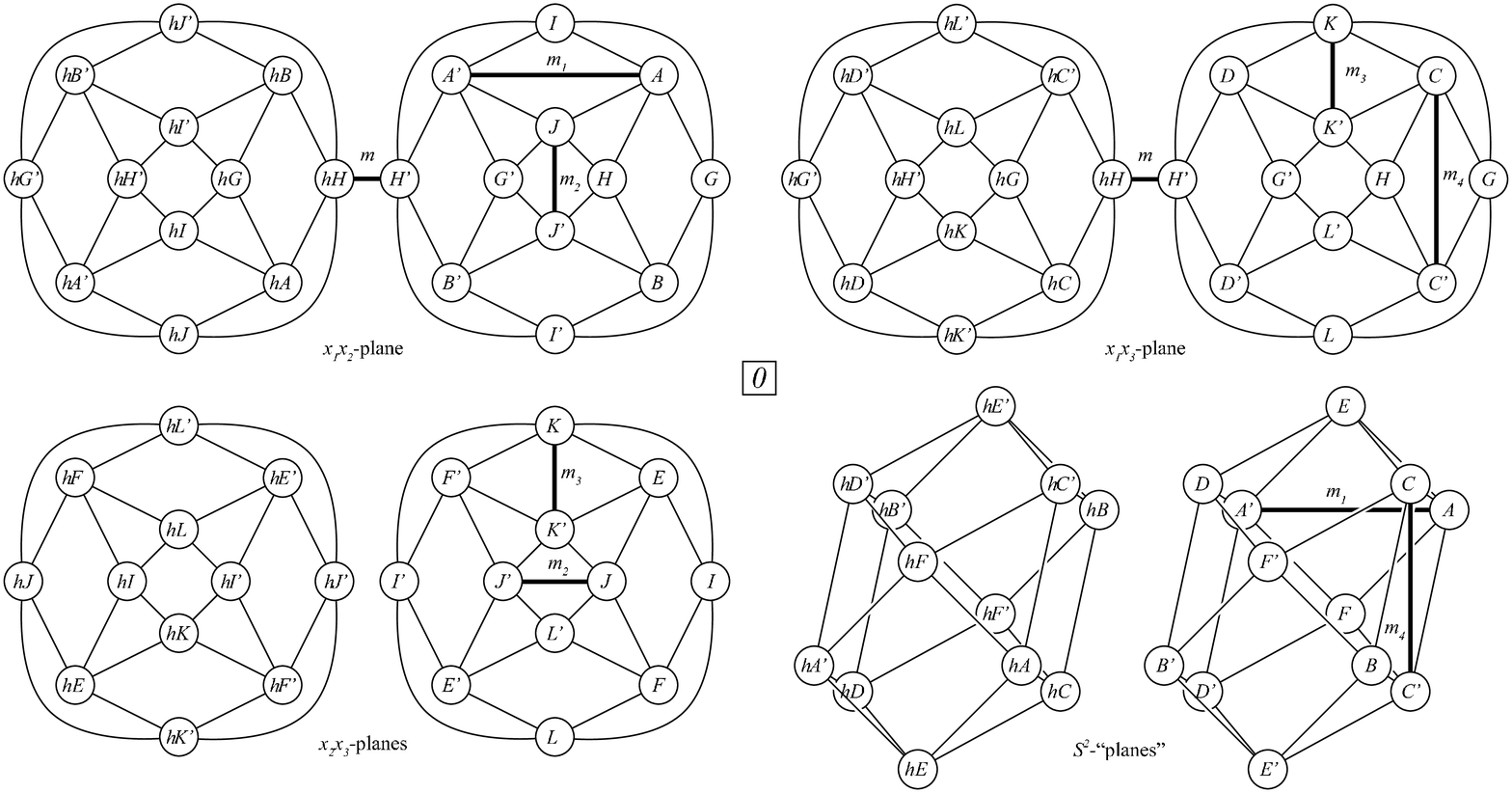}}
\caption{Step 0, initial handle decomposition}
\label{step0}
\end{figure}

This is the initial setup.
Altogether, there are twenty-four 1-handles, fifty-four 2-handles, and also thirty-four 3-handles and
five 4-handles, whose attaching spheres we do not need to keep track of.

Each of the 2-handles $m_1$, $m_2$, $m_3$ and $m_4$ (coming from attaching $T_2\times D^2$
to $\tilde M_{1011}$) passes exactly once over 1-handles $AA'$, $JJ'$, $KK'$ and $CC'$, respectively,
so those 2-handles cancel the 1-handles.
\vfill

\pagebreak

\begin{figure}
\resizebox{8in}{!}{\includegraphics{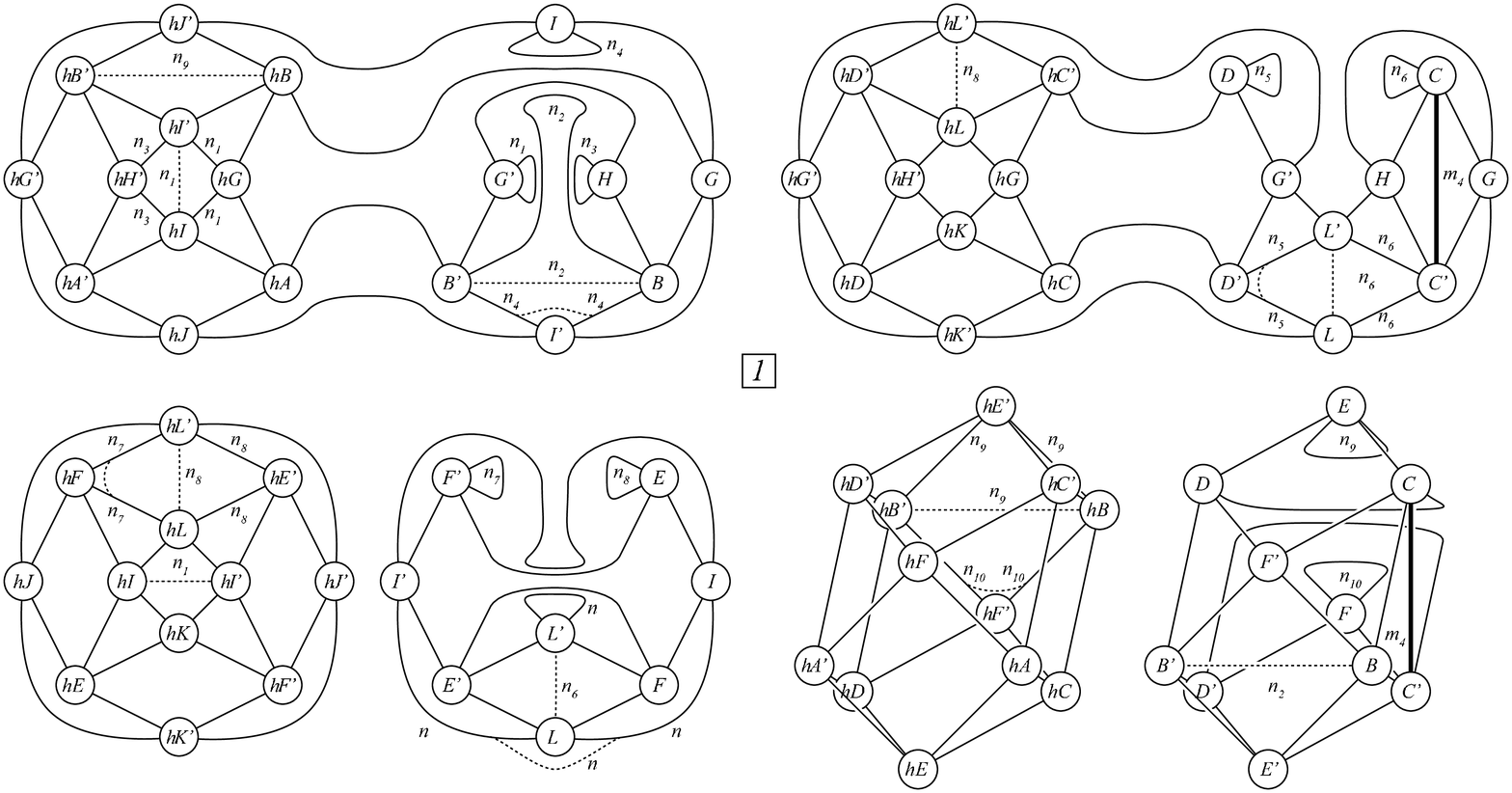}}
\caption{Step 1}
\label{step1}

\begin{tabular}{|l||l|}
\hline
\multicolumn{2}{|l|} {2-handles $m$, $m_1$, $m_2$ and $m_3$
cancel $H'hH$, $AA'$, $JJ'$ and $KK'$, respectively}
\\ \hhline{:=t::=:} 
\parbox[t]{4in}{$n_1$ slides over the 1-handle $hGG' $ and off the foot $hG$
into position indicated by the dotted line\\
$n_2$ isotopes into dotted position\\
$n_4$ slides over II', and off foot I' into dotted position\rule[-6pt]{0pt}{0pt}}
&
\parbox[t]{4in}{$n_6$ slides over $CC'$ and off foot $C'$ into dotted position\\
$n_5$ slides over DD' and off foot D' into dotted position}
\\ \hhline{:=::=:}  
\parbox[t]{4in}{$n_7$ slides over $F'hF$ and off foot $hF$ into dotted position\\
$n_8$ slides over $EhE'$ and off foot $hE'$ into dotted position\\
$n$ slides over $LL'$ and off foot $L$ into dotted position\rule[-6pt]{0pt}{0pt}}
&
\parbox[t]{4in}{$n_9$ slides over $EhE'$ and off foot $hE'$ into dotted position\\
$n_{10}$ slides over  $FhF'$ and off foot $hF'$ into dotted position}
\\ \hhline{|-||-|} 
\end{tabular}
\end{figure}

\pagebreak

\begin{figure}
\resizebox{8in}{!}{\includegraphics{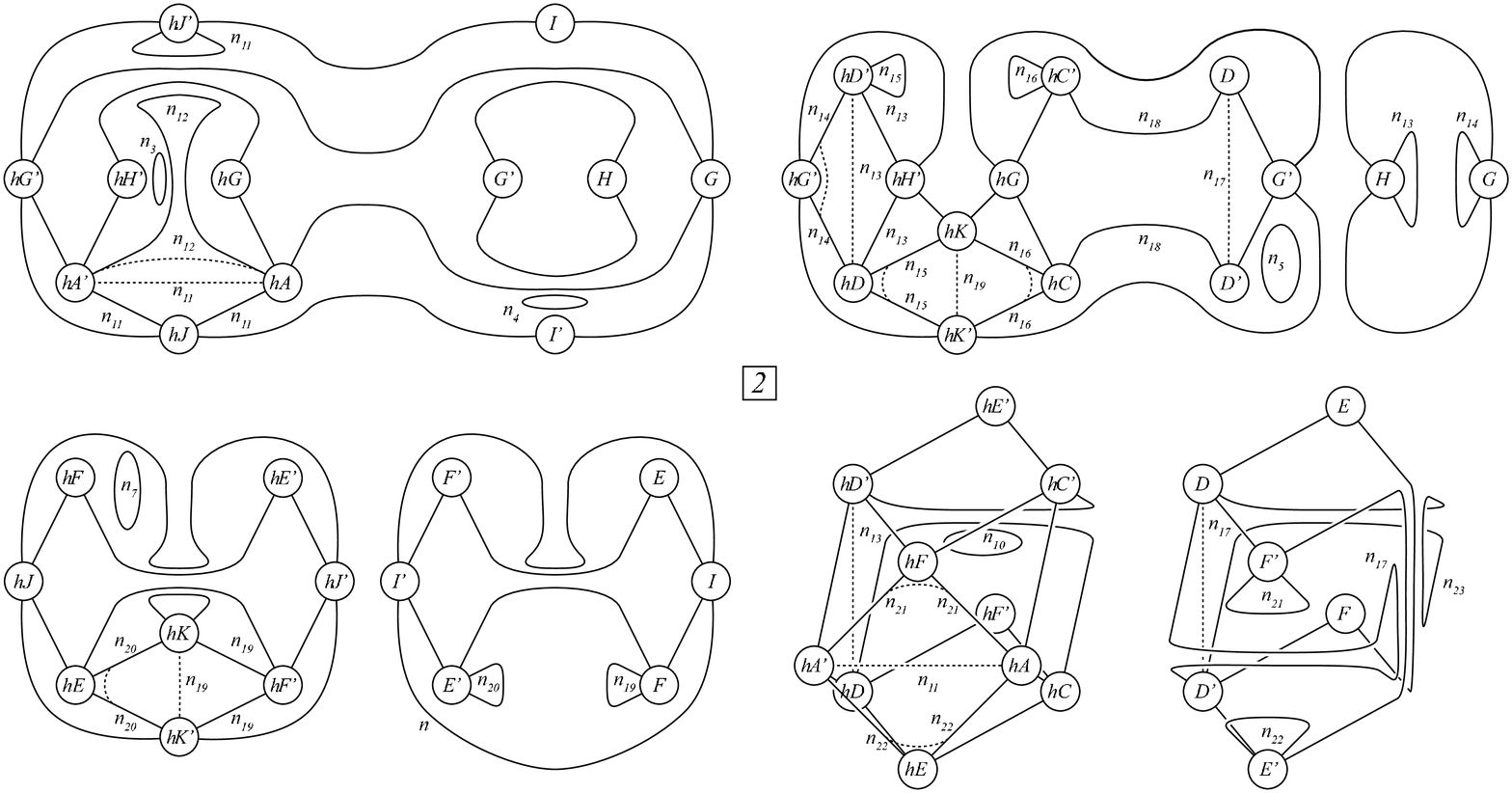}}
\caption{Step 2}
\label{step2}

\begin{tabular}{|l||l|}
\hline
\multicolumn{2}{|l|}{$m_4$, $n_1$, $n_2$, $n_6$, $n_8$ and $n_9$ cancel
$CC'$, $hIhI'$, $BB'$, $LL'$ and $hLhL'$, hBhB' respectively}
\\ \hhline{:=t::=:} 
\parbox[t]{4in}{$n_{11}$ slides over $hJhJ'$ and off $hJ$ into dotted position\\
$n_{12}$ isotopes into dotted position}
&
\parbox[t]{4in}{$n_{13}$ slides over $HhH'$ and off $hH'$ into dotted position\\
$n_{14}$ slides over $GhG'$ and off $hG'$ into dotted position\\
$n_{15}$ slides over $hDhD'$ and off $hD$  into dotted position\\
$n_{16}$ slides over $hChC'$ and off $hC'$ into dotted position\rule[-6pt]{0pt}{0pt}}
\\ \hhline{:=::=:}  
\parbox[t]{4in}{$n_{19}$ slides over $FhF'$ and off $hF'$ into dotted position\\
$n_{20}$ slides over $E'hE$ and off $hE$ into dotted position}
&
\parbox[t]{4in}{$n_{17}$ isotopes into dotted position\\
$n_{21}$ slides over $F'hF$ and off $hF$ into dotted position\\
$n_{22}$ slides over $E'hE$ and off $hE$ into dotted position\rule[-6pt]{0pt}{0pt}}
\\ \hhline{:=b::=:}
\multicolumn{2}{|l|}{
Each of the 2-handles $n_3$, $n_4$, $n_5$, $n_7$ and $n_{10}$ cancels a 3-handle owing to
Proposition~\ref{cancel3handle}}
\\ \hline
\end{tabular}
\end{figure}

\pagebreak

\begin{figure}
\resizebox{8in}{!}{\includegraphics{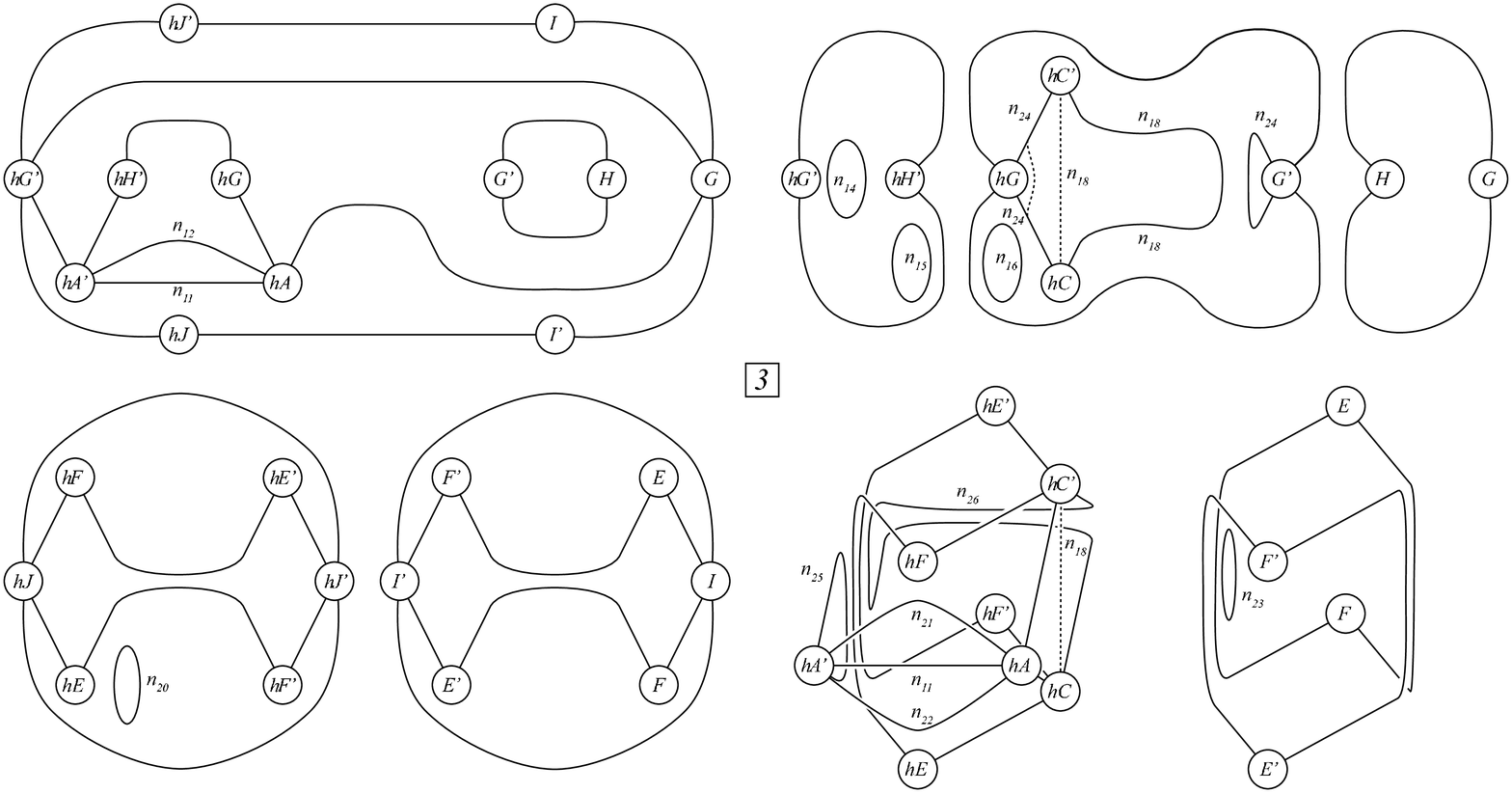}}
\caption{Step 3}
\label{step3}

\begin{tabular}{|l||l|}
\hline
\multicolumn{2}{|l|} {$n_{13}$, $n_{17}$, $n_{19}$ cancel 
$hDhD'$, $hKhK"$,  $DD'$, respectively}
\\ \hhline{:=t::=:} 
\parbox[t]{4in}{isotopy simplifies picture}
&
\parbox[t]{4in}{$n_{18}$ isotopes into dotted position\\
$n_{24}$ slides over $G'hG$ and off $hG$ into dotted position\rule[-6pt]{0pt}{0pt}}
\\ \hhline{:=::=:}  
\parbox[t]{4in}{isotopy simplifies picture}
&
\parbox[t]{4in}{$n_{23}$ isotopes\rule[-6pt]{0pt}{0pt}}
\\ \hhline{:=b::=:}
\multicolumn{2}{|l|} {$n_{14}$, $n_{15}$, $n_{16}$, $n_{20}$, and $n_{23}$
cancel a 3-handle}
\\ \hline
\end{tabular}
\end{figure}

\pagebreak

\begin{figure}
\resizebox{8in}{!}{\includegraphics{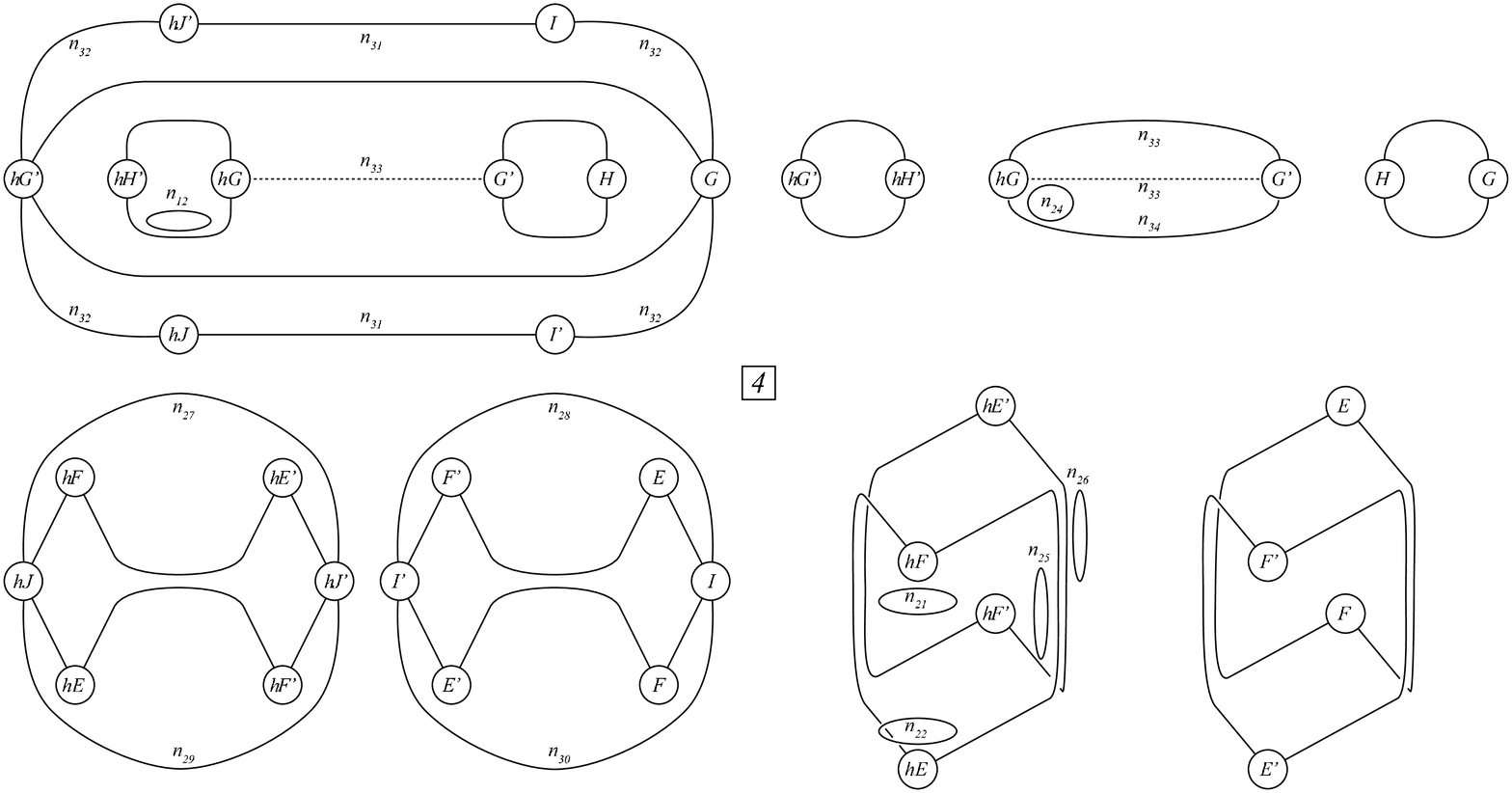}}
\caption{Step 4}
\label{step4}

\begin{tabular}{|l||l|}
\hline
\multicolumn{2}{|l|} {$n_{11}$, $n_{18}$ cancel  $hAhA'$,  $hChC'$ respectively}
\\ \hhline{:=t::=:} 
\parbox[t]{4in}{\ }
&
\parbox[t]{4in}{isotopy simplifies picture\\
$n_{33}$ isotopes to dotted position\rule[-6pt]{0pt}{0pt}}
\\ \hhline{:=::=:}  
\parbox[t]{4in}{}
&
\parbox[t]{4in}{$n_{25}$ isotopes\\
$n_{26}$ isotopes\rule[-6pt]{0pt}{0pt}}
\\ \hhline{:=b::=:}
\multicolumn{2}{|l|} {$n_{12}$, $n_{24}$, $n_{21}$, $n_{22}$, $n_{25}$, $n_{26}$ cancel a 3-handle}
\\ \hline
\end{tabular}
\end{figure}

\pagebreak

\begin{figure}
\resizebox{8in}{!}{\includegraphics{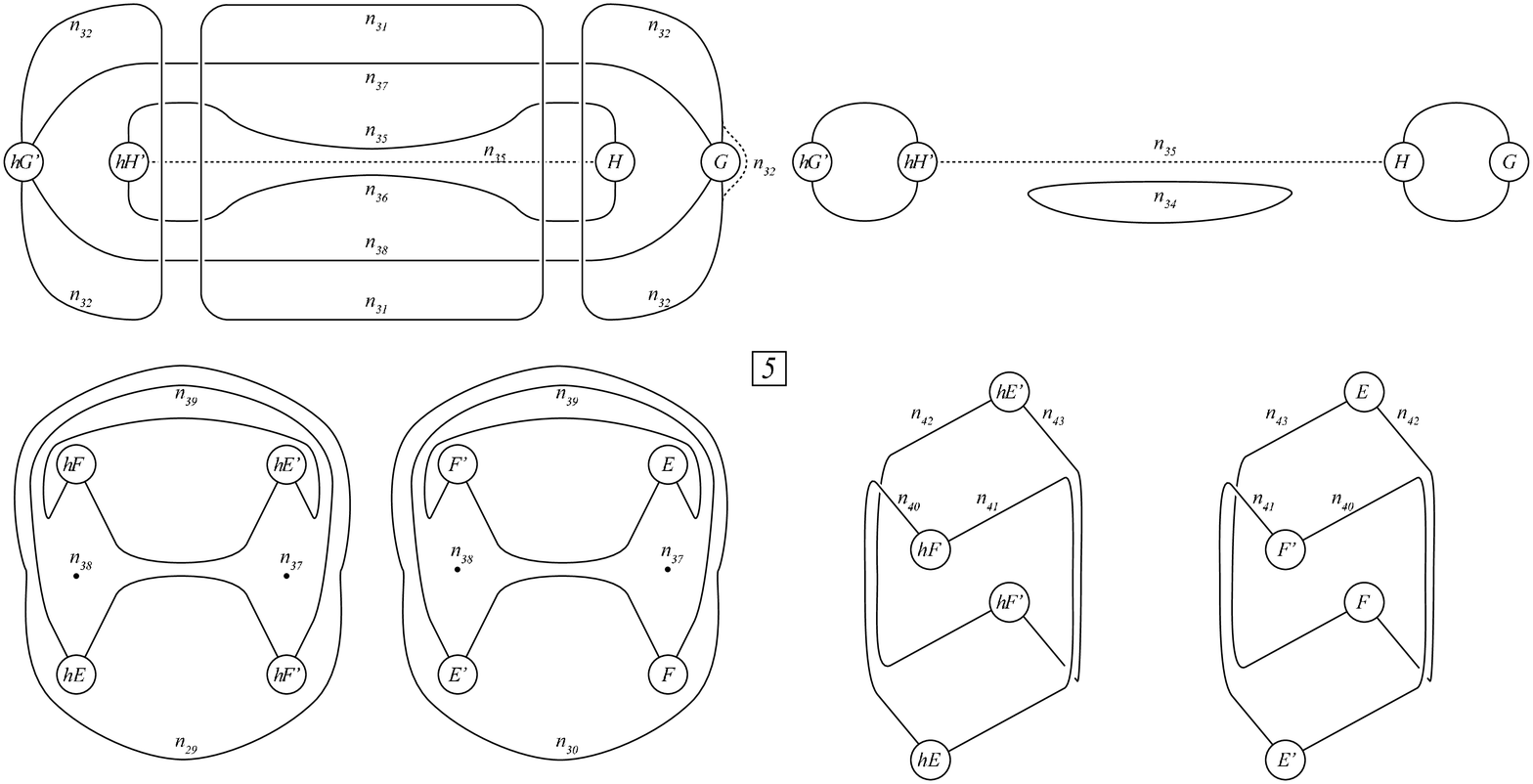}}
\caption{Step 5}
\label{step5}

\begin{tabular}{|l||l|}
\hline
\multicolumn{2}{|l|} {$n_{33}$, $n_{27}$, $n_{28}$  cancel  $G'hG$,  $hJhJ'$, $II'$ respectively}
\\ \hhline{:=t::=:} 
\parbox[t]{7in}{$n_{31}$ rises in part above $x_1x_2$-plane, spans surface like an arched roof
that contains its parallel circle\\
$n_{32}$ rises in part above $x_1x_2$-plane, slides over $GhG'$ and off $G$ into dotted position\\
$n_{35}$ isotopes into dotted position\rule[-6pt]{0pt}{0pt}}
&
\parbox[t]{1in}{}
\\ \hhline{:=::=:}  
\parbox[t]{7in}{both branches of $n_{39}$ are isotoped by rotating by $180^\circ$ around the axis joining
their endpoints --- the dots representing where $n_{37}$ and $n_{38}$ cross these planes show
they do not interfere with the isotopy\rule[-6pt]{0pt}{0pt}}
&
\parbox[t]{1in}{}
\\ \hhline{:=b::=:}
\multicolumn{2}{|l|} {\parbox[t]{8in}{$n_{31}$, $n_{32}$, $n_{34}$ cancel a 3-handle;
$n_{29}$ and $n_{30}$  can be separated from the rest of the diagram by pulling them in
 the $x_1$-direction, then each cancels a 3-handle}}
\\ \hline
\end{tabular}
\end{figure}

\pagebreak

\begin{figure}
\resizebox{8in}{!}{\includegraphics{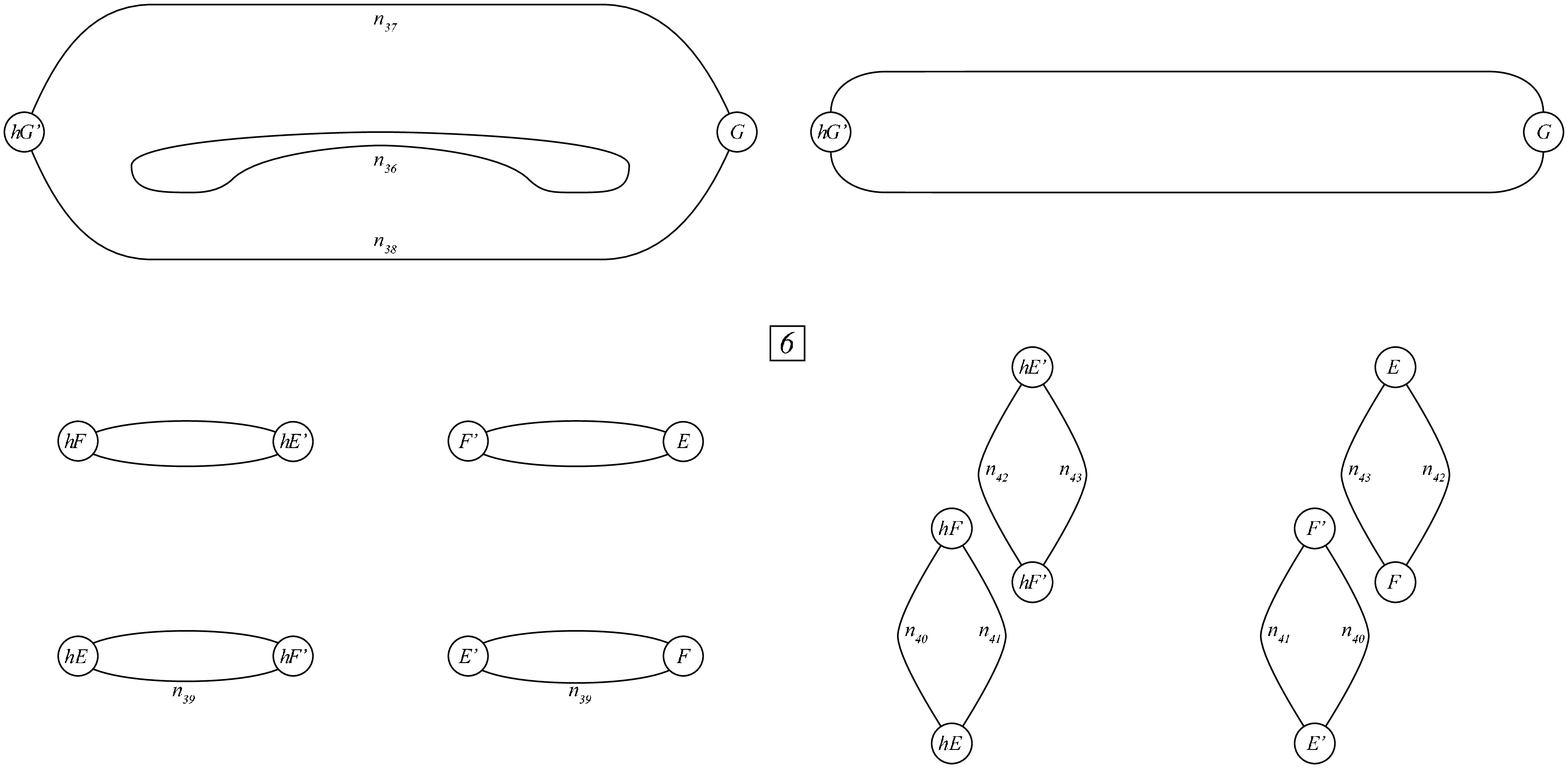}}
\caption{Step 6}
\label{step6}

\begin{tabular}{|l||l|}
\hline
\multicolumn{2}{|l|} {$n_{35}$ cancels $HhH'$}
\\ \hhline{:=t::=:} 
\parbox[t]{4in}{$n_{37}$ and $n_{38}$ isotoped up and down, respectively\rule[-6pt]{0pt}{0pt}}
&
\parbox[t]{4in}{}
\\ \hhline{:=::=:}  
\parbox[t]{4in}{}
&
\parbox[t]{4in}{$n_{40}$, $n_{41}$, $n_{42}$ and  $n_{43}$ are isotoped along the $S^2$'s so
they lie, along with their parallel circles, in planes parallel to the $x_1x_3$-plane --- note that
$n_{37}$ and $n_{38}$ do not interfere with this in their new position\rule[-6pt]{0pt}{0pt}}
\\ \hhline{:=b::=:}
\multicolumn{2}{|l|} {$n_{36}$ cancels a 3-handle}
\\ \hline
\end{tabular}
\end{figure}

\pagebreak

\begin{figure}
\resizebox{8in}{!}{\includegraphics{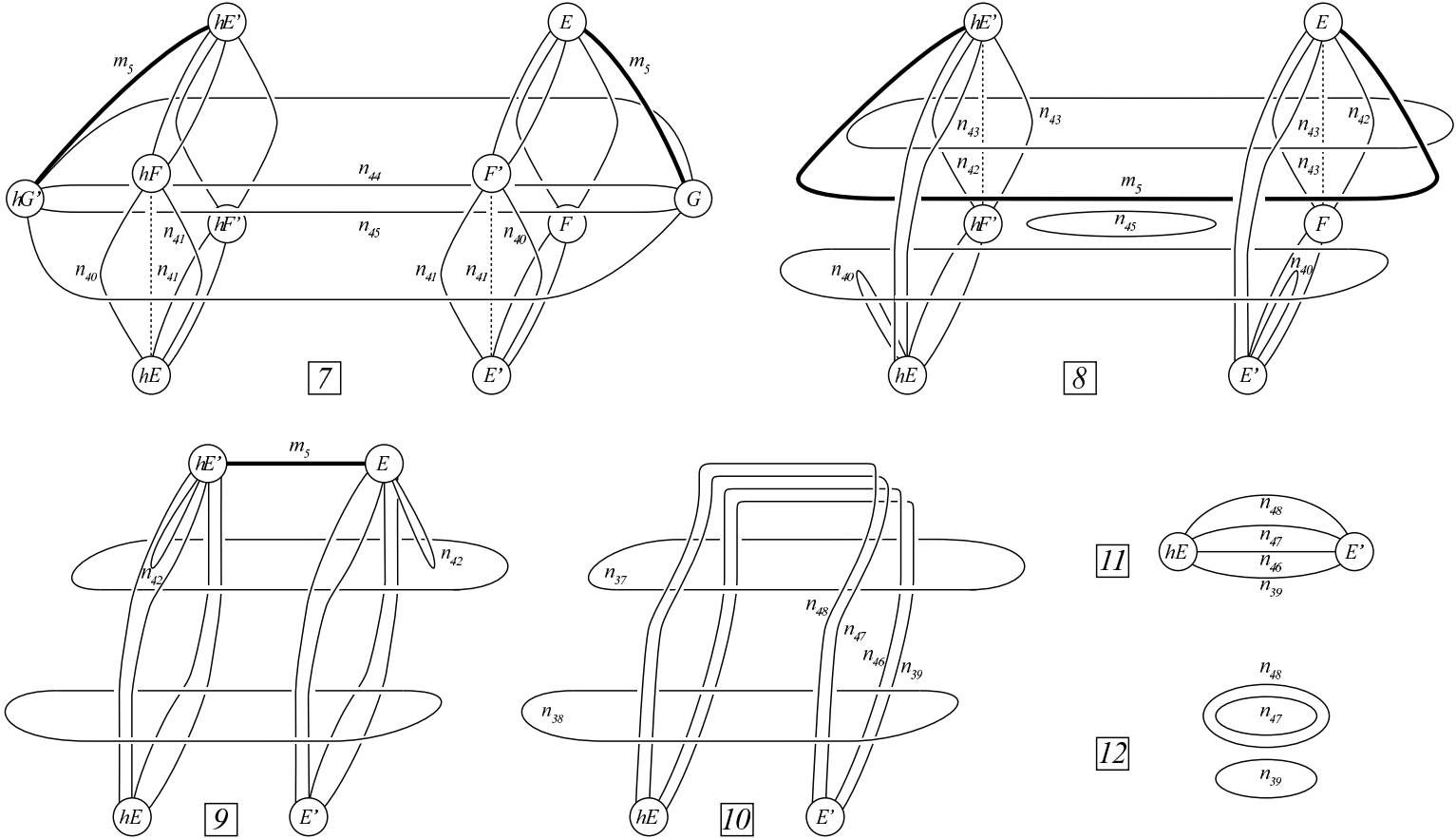}}
\caption{Steps 7--12}
\label{step7}

\begin{tabular}{|l||l|} 
\hhline{|-||-|}
\parbox[t]{4.1in}{$n_{41}$ isotopes to dotted positions\\  
parallel curve of $m_5$ can be chosen right above $m_5$}
&
\parbox[t]{4.1in}{$n_{41}$ cancels $F'hF$, $n_{44}$ cancels $GhG'$\\  
$n_{43}$ isotopes to dotted positions\\
$n_{40}$ slides over $E'hE$ and off $hE$, then cancels 3-handle\\
$n_{45}$ cancels 3-handle\rule[-6pt]{0pt}{0pt}}
\\ \hhline{|-||-|}
\end{tabular}
\begin{tabular}{|l||l||l|}
\hhline{|-||-||-|}
\parbox[t]{2in}{$m_5$ isotopes\\  
$n_{43}$ cancels $FhF'$\\
$n_{42}$ slides over $EhE'$ and off~$E$, then cancels 3-handle}
&
\parbox[t]{2.5in}{$m_5$ cancels $EhE'$\\  
$n_{37}$ and $n_{38}$ can be separated from\\
the diagram and cancel 3-handles}
&
\parbox[t]{3.5in}{$n_{39}$, $n_{46}$, $n_{47}$ and $n_{48}$ can be isotoped so they all
lie in a plane, along with their parallel curves\\  
$n_{46}$ cancels $E'hE$\\
$n_{39}$, $n_{47}$ and $n_{48}$ cancel 3-handles\rule[-6pt]{0pt}{0pt}}
\\ \hhline{|-||-||-|}
\end{tabular}
\end{figure}

\end{landscape}

\pagebreak
\bibliography{handleconv}
\bibliographystyle{plain}

\end{document}